\pdfoutput=1

\documentclass[twocolumn]{autart}    

\usepackage[utf8]{inputenc}
\usepackage{amsmath,bm}
\usepackage{amsfonts}
\usepackage{amssymb}
\usepackage{graphicx}
\usepackage{amsfonts}
\usepackage{hyperref}
\usepackage{caption}
\usepackage{cite}
\usepackage{subcaption}
\usepackage{algorithm}
\usepackage[noend]{algpseudocode}
\usepackage{xcolor}
\usepackage{enumerate}
\usepackage[normalem]{ulem}

\newtheorem{definition}{Definition}
\newtheorem{Theorem}{Theorem}

\newtheorem{Lemma}{Lemma}
\newtheorem{Problem}{Problem}
\newtheorem{Corollary}{Corollary}
\newtheorem{Remark}{Remark}
\newtheorem{Assumption}{Assumption}

\usepackage{tikz}
	\usetikzlibrary{shapes,arrows}
	\tikzstyle{frame} = [draw, -latex]
	\tikzstyle{line} = [draw]
	\tikzstyle{line2} = [draw, dashdotted]
	\tikzstyle{line3} = [draw, dashed]
	\tikzstyle{line3UD} = [draw, dashed]
	\tikzstyle{place} = [circle, draw=black, fill=white, thick, inner sep=2pt, minimum size=1mm]
	\tikzstyle{place2} = [circle, draw=black, fill=black, thick, inner sep=2pt, minimum size=1mm]
	\tikzstyle{placeRed} = [circle, draw=red, fill=red, thick, inner sep=2pt, minimum size=1mm]
	\tikzstyle{vertex} = [circle, draw=black, fill=black, thick, inner sep=2pt, minimum size=1mm]
\usepackage{tkz-euclide}

\newcommand{\m}[1]{\mathbf{#1}}
\newcommand{\mc}[1]{\mathcal{#1}}
\newcommand{\mb}[1]{\mathbb{#1}}
\newcommand{\abs}[1]{\lVert{#1} \rVert}
\newcommand\numberthis{\addtocounter{equation}{1}\tag{\theequation}}
\newcommand*{\eprf}{\hfill\ensuremath{\square}}
\makeatletter
\def\algbackskip{\hskip-\ALG@thistlm}
\def\@opargbegintheorem#1#2#3{\trivlist
   \item[]{\bfseries #1\ #2\ (#3)} \itshape}
\makeatother
\usepackage[switch,pagewise]{lineno}
\newcommand*\patchAmsMathEnvironmentForLineno[1]{%
\expandafter\let\csname old#1\expandafter\endcsname\csname #1\endcsname
\expandafter\let\csname oldend#1\expandafter\endcsname\csname end#1\endcsname
\renewenvironment{#1}%
{\linenomath\csname old#1\endcsname}%
{\csname oldend#1\endcsname\endlinenomath}}%
\newcommand*\patchBothAmsMathEnvironmentsForLineno[1]{%
\patchAmsMathEnvironmentForLineno{#1}%
\patchAmsMathEnvironmentForLineno{#1*}}%
\AtBeginDocument{%
\patchBothAmsMathEnvironmentsForLineno{equation}%
\patchBothAmsMathEnvironmentsForLineno{align}%
\patchBothAmsMathEnvironmentsForLineno{flalign}%
\patchBothAmsMathEnvironmentsForLineno{alignat}%
\patchBothAmsMathEnvironmentsForLineno{gather}%
\patchBothAmsMathEnvironmentsForLineno{multline}%
}

\begin{document}

\begin{frontmatter}

\title{Pose Localization of Leader-Follower Networks with Direction Measurements\thanksref{footnoteinfo}} 
\thanks[footnoteinfo]{This paper was not presented at any IFAC 
meeting. Corresponding author Hyo-Sung Ahn. Tel. +82-62-970-2398.}
\author[GIST]{Quoc Van Tran}\ead{tranvanquoc@gist.ac.kr},    
\author[ANU]{Brian D. O. Anderson}\ead{Brian.Anderson@anu.edu.au},               
\author[GIST]{Hyo-Sung Ahn}\ead{hyosung@gist.ac.kr}  
\address[GIST]{School of Mechanical Engineering,
        Gwangju Institute of Science and Technology, Gwangju, Republic of Korea.}  
\address[ANU]{Hangzhou Dianzi University, Hangzhou 310018, China, Data61-CSIRO,
and Research School of Electrical, Energy and Materials Engineering, The Australian National University,
Canberra ACT 2601, Australia.}             


\begin{abstract}                          
A distributed pose localization framework based on direction measurements is proposed for a type of \textit{leader-follower} multi-agent systems in $\mb{R}^3$. 
The novelty of the proposed localization method lies in the elimination of the need for using distance measurements and relative orientation measurements for the network pose localization problem.
In particular, a network localization scheme is developed based directly on the measured direction constraints between an agent and its neighboring agents in the network. The proposed position and orientation localization algorithms are implemented through differential equations which simultaneously compute poses of all followers by using locally measured directional vectors and angular velocities, and actual pose knowledge of some leader agents, allowing some tracking of time-varying orientations. Further, we establish an almost global asymptotic convergence of the estimated positions and orientations of the agents to the actual poses in the stationary case. 
\end{abstract}

\end{frontmatter}

\section{Introduction}
The network localization is a crucial step that often needs to be done first in order for multi-agent systems to perform further coordination control or distributed estimation tasks \cite{Oh2015survey,Zhao2018csm}. Distances and directions are the two most commonly used measurements that are widely used in position localization literature \cite{Aspnes2006, Mao2007}. In a three dimensional ambient space, direction is characterized by a unit length vector; thus the directional vector to a target can be easily computed from its pixel coordinates in an image by using the pin-hole camera model \cite{Ma2004} when a visual image is available. Furthermore, in three-dimensional space, additional relative orientation\footnotemark{} measurements between neighboring agents are often required for estimating orientations of the agents in a network, a process which is called orientation localization \cite{Piovan2013,Tron2014tac,Montijano2016tro}. \footnotetext{A relative orientation is effectively the rotation matrix linking a local coordinate frame of one agent to the local coordinate frame of another agent. It is often estimated by vision-based techniques, e.g. by processing images (of a common scene) captured by the agents and establishing the feature correspondences \cite{Tron2014tac}.} However, there are not many works that study simultaneous localization of positions and orientations, a process which is called pose localization, in a distributed setup. 
Motivated by these facts, this work attempts to provide a distributed pose localization framework for a type of leader-follower networks. Moreover, in order not to use the relative orientation measurements, this paper uses the direction measurements and pose knowledge of some leader agents.

For a two-dimensional ($2$-D) ambient space,  network localization laws using angles of arrival between triplets of nodes are proposed in \cite{Zhu2008} and an orientation localization method utilizing orientation knowledge of a few nodes is presented in \cite{Rong2006}. The authors in \cite{Piovan2013} further proposed a least-squares optimization problem to achieve orientation localization by exploiting kinematic relationships among the orientations of nodes. A least-squared algorithm for position localization using bearing-only information is proposed in \cite{Bishop2011ifac}. In $3$-dimensional space ($3$-D), it is often required that relative orientation measurements are available for estimating the orientations of the agents. For example, some necessary and sufficient conditions are provided for orientation localizability of triangular sensing networks of relative orientation measurements in \cite{Piovan2013}, without providing a distributed orientation localization law. Network localization schemes using relative poses (relative orientations and relative positions), which are measured by a vision-based technique, are investigated in \cite{Tron2014tac,Thunberg2017cdc}. The estimation of relative poses, however, requires the agents to have views of a common scene and complicated estimation algorithms. By using relative orientation measurements, our recent works in \cite{Quoc2018cdc, Quoc2018tcns} propose distributed orientation estimation laws which guarantee almost global convergence of the estimated orientations up to a common orientation. The authors in \cite{Zhao2016aut} propose a direction-only position localization law for \textit{bearing rigid} networks with two anchor nodes. However, \cite{Zhao2016aut} further assumes that the agents know their actual orientations. There is no framework for direction-only network localization and formation control in $3$-D when agents lack knowledge of a global frame.

The orientation localization problem is challenging and requires sophisticated estimation algorithms.
In $2$-D, it is straightforward to see how two neighboring agents observing each other might determine a common view of their relative orientation (i.e., a scalar angle), within an unknown constant rotation common to both, see e.g. \cite{Lee2016auto,Oh2014tac, Wang2018tsp}, as is now described. Each agent maintains a (possibly body-fixed) coordinate frame and measures the orientation angle of its neighboring agent (assuming direction sensing technology). In any common frame, the measured angles (of the two neighboring agents) must differ by precisely $\pi$ radians. Hence a rotation of the coordinate axes of one agent can be made to ensure that after rotation, the angle difference is compensated. For an $n$ agent network, one has to put together in a distributed fashion a collection of such calculations.

\begin{figure}[t]
\centering
\begin{tikzpicture}
\node (pz) at (0,1,0) []{};
\node (bij) at (1.15,0.5,0) {};
\node (bji) at (1.2,.5,0) {};
\node (p_j) at (2.5,1,0) [label=below:$\m{p}_j$]{};
\node (w_i) at (-0.4,.6,0) [label=left:${^i\Sigma}$]{};
\node (Sigj) at (2.4,1.7,0) [label=right:${^j\Sigma}$]{};

\draw[{line width=0.7pt},blue] (0,0,0) [frame] -- (0.7,0,0);
\draw[{line width=0.7pt},blue,->] (0,0,0)[frame]   -- (w_i);
\draw[{line width=0.7pt},blue] (0,0,0) [frame] -- (pz);
\draw[{line width=0.7pt},blue] (p_j) [frame] -- (1.95,1.4,0);
\draw[{line width=0.7pt},blue,->] (p_j)[frame]  -- (3,1.4,0);
\draw[{line width=0.7pt},blue] (p_j) [frame] -- (2.5,1.8,0);

\draw[{line width=1.5pt}] (0,0,0) [frame] -- (bij) node [pos =0.9, yshift=-0.2ex, above left] {$\m{b}_{ij}^i$};
\draw[{line width=1.5pt}] (p_j) [frame] -- (bji) node [pos =0.75, yshift=-0.6ex, above left] {$\m{b}_{ji}^j$};
\draw[{line width=.7pt},black] (0,0,0)  -- (p_j);

\node[place] (p_i) at (0,0,0) [label=left:$\m{p}_i$]{};
\node[place] () at (p_j) []{};
\end{tikzpicture}
\caption{The agents $i$ and $j$ respectively measure the directions $\m{b}_{ij}^i\in \mb{R}^3$ and $\m{b}_{ji}^j \in\mb{R}^3$ in local coordinate frames. Using these measurements, they would like to decide the relative orientation $\m{R}_{ij}\in SO(3)$ and the orientations $\m{R}_i$ and $\m{R}_j$.} 
\label{fig:relative_orient_goal}
\end{figure}
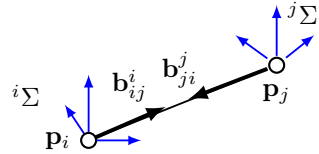

How to do something like this in a $3$-dimensional ambient space is less clear. For example, with only a pair of direction measurements between two neighboring agents $i$ and $j$ ($\m{b}_{ij}^i,\m{b}_{ji}^j$) $\in \mb{R}^3\times \mb{R}^3$ (see Fig. \ref{fig:relative_orient_goal}), it is insufficient for the agents $i$ and $j$ to determine their \textit{relative orientation}, i.e., $\m{R}_{ij}\triangleq \m{R}_i^\top\m{R}_j\in SO(3)$, where $\m{R}_i$ and $\m{R}_j\in SO(3)$ are the orientation matrices of agents $i$ and $j$, respectively, due to the ambiguity of the rotation along the common direction vector, $\m{b}_{ij}$. This difficulty can be overcome by examining additional direction constraints of each of the two agents to a third agent $k$ that they both observe. Indeed, as shown in \cite{Quoc2018cdc}, by exploiting the triangle sensing network and using a coordinate frame alignment procedure, agents $i$ and $j$ can compute $\m{R}_{ij}$. The orientations of all agents then can be computed up to a common orientation bias by using a consensus protocol \cite{Quoc2018cdc}. This method, however, relies on the existence of triangle networks and requires predefinition of a complicated computation sequence. 

This paper proposes a distributed  pose localization scheme for a type of leader-follower network that uses continuous-time direction vectors and two or more anchor agents which know their absolute poses.
A distributed orientation localization protocol in $SO(3)$ that estimates orientations of all followers is proposed. Under the proposed orientation localization protocol, estimated orientations converge to the true orientations of agents almost globally and asymptotically. By using the estimates of orientations and direction measurements, we investigate a position localization law for the \textit{leader-follower} network. Under the proposed position localization law, positions of all followers are also globally and asymptotically determined. The proposed network pose localization laws can work exclusively with inter-agent directional vectors and does not require a common scene and extra algorithms to compute relative poses; unlike \cite{Tron2014tac,Thunberg2017cdc}.

The rest of this paper is organized as follows. Section \ref{sec:preliminary} presents some preliminaries and the problem formulation. The orientation localization problem is studied in Section \ref{sec:orient_local}. We proposes a position localization law and establish the global asymptotic convergence of computed positions in Section \ref{sec:pos_local}. 
Finally, numerical simulations are provided in Section \ref{sec:Sim} and Section \ref{sec:Conclusion} concludes this paper.

\section{Preliminaries and Problem Formulation}\label{sec:preliminary} 
In this paper we use the following notations. 
The dot product and cross product are denoted by $\cdot$ and $\times$, respectively. The symbol $\Sigma$ represents a global coordinate frame and the symbol $^k\Sigma$ with the superscript index $k$ denotes the $k$-th local coordinate frame. Let $\m{1}_n=[1,\ldots,1]^\top\in \mb{R}^n$ be the vector of all ones, and $\m{I}_3$ the $3\times 3$ identity matrix. The Kronecker product is denoted by $\otimes$. The trace of a matrix is denoted by $\text{tr}(\cdot)$. The set of rotation matrices in $\mb{R}^3$ is denoted by $SO(3)=\{\m{Q}\in\mb{R}^{3\times 3}~|~\m{Q}\m{Q}^\top=\m{I}_3,\text{det}(\m{Q})=1\}$. The set of real matrices with orthonormal column vectors is $O(3)$. The orthogonal projection matrix associated with a nonzero vector $\m{x}\in \mb{R}^{3}$ is defined as
\begin{equation}\label{eq:projection_matrix}
\m{P}_{\m{x}} = \m{I}_3-\frac{\m{x}}{||\m{x}||}\frac{\m{x}^\top}{||\m{x}||}\in \mb{R}^{3\times 3}.
\end{equation}
It can be verified that $\m{P}_{\m{x}}$ is positive semidefinite and idempotent. Moreover, $\m{P}_{\m{x}}$ has the nullspace $\text{null}(\m{P}_{\m{x}})=\text{span}\{\m{x}\}$ and the eigenvalue set $\{0,1,1\}$ \cite{Zhao2018csm}.

The space of $3\times 3$ skew-symmetric matrices is denoted by $\mathfrak{so}(3):=\{\m{A}\in \mb{R}^{3\times 3}|\m{A}^\top=-\m{A}\}$. For any $\omega=[\omega_x,\omega_y,\omega_z]^\top\in \mb{R}^3$, the \textit{hat} map $(\cdot)^{\wedge}:~\mb{R}^3\rightarrow \mathfrak{so}(3)$ is defined such that $\omega\times \m{v}=\omega^\wedge\m{v},\forall \m{v}\in \mb{R}^3$.
The \textit{vee} map is the inverse of the \textit{hat} map and defined as $(\cdot)^\vee:~\mathfrak{so}(3)\rightarrow \mb{R}^3$ \cite{Bullo2005spr}. The \textit{exponential map} $exp:\mathfrak{so}(3)\rightarrow SO(3)$ is \textit{surjective} and $T_{\m{R}}SO(3)=\{\m{R}\eta^\wedge:\eta^\wedge\in \mathfrak{so}(3)\}$ denotes the tangent space at a point $\m{R}\in SO(3)$.

For any $\m{x},\m{y}\in \mb{R}^3$, $\m{A},\m{B}\in \mb{R}^{3\times 3}$, and $\m{R}\in SO(3)$ we have \cite{Bullo2005spr, Mahony2008tac}
\begin{align}
&\qquad\m{x}\times \m{y}=-\m{y}\times \m{x}  \label{eq:cross_prod_1} \\
&(\m{R}\m{x})\times(\m{R}\m{y})=\m{R}(\m{x}\times\m{y}),~\m{R}\m{x}^\wedge\m{R}^\top=[\m{Rx}]^\wedge \label{eq:cross_prod_2}\\
&\qquad (\m{x}\times\m{y})^\wedge=\m{x}^\wedge\m{y}^\wedge-\m{y}^\wedge\m{x}^\wedge=\m{y}\m{x}^\top-\m{x}\m{y}^\top\label{eq:cross_prod_3}\\
&\qquad\m{x}^\top\m{y} = \text{tr}(\m{x}\m{y}^\top) \label{eq:trace_property_1}\\
&\text{tr}(\m{x}^\wedge\m{A}) = \frac{1}{2}\text{tr}[\m{x}^\wedge(\m{A}-\m{A}^\top)]=-\m{x}^\top(\m{A}-\m{A}^\top)^\vee\label{eq:trace_property_2} \\
&\qquad \text{tr}(\m{A}\m{B})=\text{tr}(\m{B}\m{A})=\text{tr}(\m{A}^\top\m{B}^\top)\label{eq:trace_property_3}
\end{align}

\subsection{Directional vector and orientation of agent}\label{subsec:orientation_def}
Consider a network of $n$ nodes in $3$-dimensional space. Each node corresponds to an agent, and an agent is defined by the position of its centroid and the orientation of a body-fixed coordinate frame $^i\Sigma$ relative to a global frame $\Sigma$. In the sequel, the position of an agent will be taken to be the position of its centroid. 
Let $\m{p}_i$ and $\m{p}_i^i \in \mb{R}^3$ be the position of agent $i$ expressed in the global frame $\Sigma$ and its body-fixed coordinate frame $^i\Sigma$, respectively. We define the unit directional vector (expressed in $\Sigma$) pointing from agent $i$ toward its neighbor $j$ along the direction of $\m{p}_{ij}$ ($\m{p}_{ij}= \m{p}_j - \m{p}_i$) as
\begin{equation*}\label{eq:bearing_vector}
\m{b}_{ij}\triangleq \frac{\m{p}_j-\m{p}_i}{\abs{\m{p}_j-\m{p}_i}}=\frac{\m{p}_{ij}}{\abs{\m{p}_{ij}}}.
\end{equation*}
The directional vector with the reverse direction is $\m{b}_{ji}=-\m{b}_{ij}$ and it points from agent $j$ toward $i$. The directional vector $\m{b}_{ij}$ measured locally in $^i\Sigma$ is denoted as $\m{b}^i_{ij}$.

The orientation or attitude of agent $i$ in $\mb{R}^3$ can be characterized by a square, orthogonal matrix $\m{R}_i\in SO(3)$ whose column vectors represent the coordinates of the orthogonal bases of the $i$-th local coordinate frame expressed in the global coordinate frame. 
The pair $(\m{R}_i,\m{p}_i) \in SE(3)$ characterizes the \textit{pose} of each agent $i$ in the global Cartesian space.

\subsection{Graph theory}
An interaction graph characterizing an interaction topology of a multi-agent network is denoted by $\mc{G}=(\mc{V},\mc{E})$, where, $\mc{V}=\{1,\ldots,n\}$ denotes the vertex set and $\mc{E}\subseteq\mc{V}\times \mc{V}$ denotes the set of edges of $\mc{G}$. An edge is defined by the ordered pair $e_k=(i,j), k=1,\ldots,m,m=\vert \mathcal{E} \vert$. The graph $\mc{G}$ is said to be undirected if $(i,j)\in \mc{E}$ implies $(j,i)\in \mc{E}$, i.e. if $j$ is a neighbor of $i$, then $i$ is also a neighbor of $j$. If the graph $\mc{G}$ is directed, $(i,j)\in \mc{E}$ does not necessarily imply $(j,i)\in \mc{E}$. The set of neighboring agents of $i$ is denoted by $\mc{N}_i=\{j\in\mc{V}:(i,j)\in \mc{E}\}$. 

\subsection{Problem formulation}
\begin{figure}[t]
\centering
\begin{subfigure}[b]{0.23\textwidth}
\centering
\begin{tikzpicture}[scale=1.]
\node[place,black] (p1) at (0,0) [label=left:$1$]{};
\node[place,black] (p2) at (1,0) [label=right:$2$]{};
\node[place] (p3) at (0.5,-0.75) [label=below right:$3$]{};
\node[place] (p4) at (-0.5,-0.75) [label=left:$4$]{};
\node[place] (p5) at (1.5,-0.75) [label=right:$5$]{};
\node[place] (p6) at (.5,-1.5) [label=left:$6$]{};
\node[place, black] (p1prime) at (1.4,-1.3) [label=right:$1'$]{};
\draw[line width=1pt,->] (p3)[frame]  -- (p1);
\draw[line width=1pt,->] (p3)[frame]  -- (p2);
\draw[line width=1pt,->] (p4)[frame]  -- (p3);
\draw[line width=1pt,->] (p4)[frame]  -- (p1);
\draw[line width=1pt,->] (p5)[frame]  -- (p1);
\draw[line width=1pt,->] (p5)[frame]  -- (p3);
\draw[line width=1pt,->] (p5)[frame]  -- (p2);
\draw[line width=1pt,->] (p4)[frame]  -- (p2);
\draw[line width=1pt,->] (p6)[frame]  -- (p4);
\draw[line width=1pt,->] (p6)[frame]  -- (p5);
\draw[line width=1pt,->] (p6)[frame]  -- (p3);
\draw[line width=1pt,->] (p6)[frame]  -- (p1prime);
\end{tikzpicture}
\caption{}
\label{fig:vertex_addition}
\end{subfigure}
\begin{subfigure}[b]{0.23\textwidth}
\centering
\begin{tikzpicture}[scale=.75]
\node (pz) at (0,1,0) []{};
\node (bij) at (0.4,0.6) [label=right:$\m{b}_{ij}^i$]{};
\node (bji) at (0.6,1.2,0) [label=right:$\m{b}_{ji}^j$]{};
\node (p_j) at (1,2,0) [label=right:$j$]{};
\draw[{line width=1.pt},red] (0,0,0) [frame] -- (0.5,1.);
\draw[{line width=1.pt},red] (p_j) [frame] -- (0.5,1.);
\draw[{line width=.5pt},black,->] (0,0,0)[frame]  -- (p_j);
\node[place] (p_i) at (0,0,0) [label=left:$i$]{};
\node[place] (p_k) at (-1,2,0) [label=left:$k$]{};
\node[place] () at (p_j) []{};
\draw[{line width=.5pt},black] (0,0,0)[frame]  -- (p_k);
\draw[{line width=1.pt},red] (0,0,0) [frame] -- (-.4,0.8,0);
\node (bik) at (-.4,.6) [label=left:$\m{b}_{ik}^i$]{};
\node (bji) at (-0.4,1.2) [label=left:$\m{b}_{ki}^k$]{};
\draw[{line width=1.pt},red] (p_k) [frame] -- (-0.5,1);
\end{tikzpicture}
\caption{}
\label{fig:bearing_constraints}
\end{subfigure}
\caption{(a) A \textit{twin-leader-follower} network constructed from the \textit{vertex addition}: anchor nodes $1$, $2$, and $1'$ (black), the first follower $3$. (b) Agent $i$ measures $(\m{b}_{ij}^i,\m{b}_{ik}^i),~j,k\in \mc{N}_i$ and receives $(\hat{\m{R}}_j\m{b}_{ji}^j,\hat{\m{R}}_k\m{b}_{ki}^k)$ from $j$ and $k$.}
\label{fig:leader_follower_def}
\end{figure}
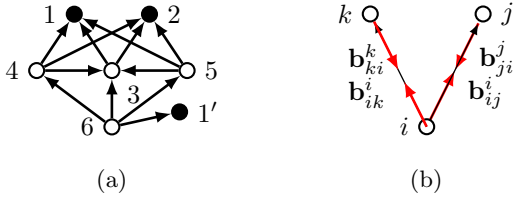

Consider a leader-follower network in $\mb{R}^3$ with at least two non-collocated leader agents $1$ and $2$ which know their actual poses (position and orientation in a global coordinate frame). 
Starting with the leader agents, the leader-follower network studied in this work is defined follows (See also Fig. \ref{fig:vertex_addition}).
\begin{definition}[Twin-Leader-Follower Network]\label{def:leader_follower}
A twin-leader follower network is a directed network in which agents are ordered such that (a) all leader agents appear first, there are two (or more) leaders $1$ and $2$ which know their absolute poses $(\m{R}_1,\m{p}_1)$ and $(\m{R}_2,\m{p}_2)$, respectively (b) a follower agent $i,~ 3\leq i \leq n,$ has at least two neighboring agents $j$'s in the set $\{1,\ldots,i-1\}$, i.e., $|\mc{N}_i|\geq 2$, where $\mc{N}_i$ denotes the set of neighboring agents of $i$. Agent $i$ knows the direction $\m{b}^i_{ij}$ to the neighbor $j$, while its neighbor knows the direction $\m{b}^j_{ji}$.
\end{definition}

We remark that the first listed nonleader agent is known as a first follower and any leader agents beyond the first two are known as redundant leaders. To streamline nomenclature, we number the agents as $\{1,2,1^\prime,2^\prime,\ldots,3,4,\ldots,n\}$, where the follower agents are $3,4,\ldots,n$; also $\mc{V}_l=\{1,2,1^\prime,2^\prime\ldots\}$, where $1^\prime,2^\prime\ldots$ are redundant leaders, and $\mc{V}_f=\{3,4,\ldots,n\}$ will denote the sets of leader and follower agents, respectively.
In Fig. \ref{fig:vertex_addition}, there is one additional redundant leader agent $1^\prime$ to which agent $6$ measures the direction. The absolute poses of the redundant leaders can be used, in addition, as reliable measurements. On the other hand, with only one leader, it is insufficient to compute the actual poses of the agents. This is due to the fact that there are translational and scale ambiguities in networks with direction-only measurements \cite{Zhao2018csm} (See also \cite[Lm. 2 $\&$ 3]{Minh2017tac}).

Each agent $i\in \mc{V}_f$ in the network aims to estimate its actual pose, i.e., $({\m{R}}_i,{\m{p}}_i)\in SO(3)\times \mb{R}^3$, based on the direction constraints to its neighboring agents and the actual poses of the leader agents. At each time instant $t$ agent $i$ holds an estimate of its pose, denoted as $(\hat{\m{R}}_i,\hat{\m{p}}_i)\in SO(3)\times \mb{R}^3$. For this, we further make an assumption on the direction measurements as follows.

\begin{Assumption}\label{ass:bidirection_sensing}
Each agent $j$ estimates its orientation at time $t$ by $\hat{\m{R}}_j$, and transmits the information $\hat{\m{R}}_j\m{b}_{ji}^j$ to agent $i,j\in\mc{N}_i$ (see Fig. \ref{fig:bearing_constraints}).
\end{Assumption}

We assume that the agents in the network do not translate but they might rotate according to the kinematics
\begin{equation*}
\dot{\m{R}}_i=\m{R}_i(\omega_i^i)^\wedge,~\text{for } i \in \mc{V},
\end{equation*}
where $\omega^i$ is the angular velocity of agent $i$ measured locally in $^i\Sigma$. We assume that $\omega^i_i$ and its derivative are bounded, i.e., $||\omega^i_i||\leq \bar{\omega}_i$, $||\dot{\omega}^i_i||\leq \dot{\bar{\omega}}_i$, for positive constants $\bar{\omega}_i,\dot{\bar{\omega}}_i>0$, and each agent $i$ can measure $\omega^i_i$ without noise. The angular velocity expressed in the global coordinates is $\omega_i=\m{R}_i\omega_i^i$. Consequently, a locally measured direction $\m{b}_{ij}^i$ is not necessarily constant due to the rotation of agent $i$, although, in global coordinates, $\m{b}_{ij}$ is constant. This kind of system might represent a visual sensor network \cite{Tron2014tac} or a system of autonomous agents in a desired formation \cite{Quoc2018tcns} where the agents might rotate to track objects. 

Moreover, for the uniqueness of the localized poses of the agents, we have the following assumption.
\begin{Assumption}\label{ass:non_colocated}
No two agents are collocated and each follower $i\in \mc{V}_f$ has at least one pair of neighbors with which it is not collinear.
\end{Assumption}

We first address the problem of calculating the orientation $\hat{\m{R}}_{i}$ for all follower agents based on the direction and angular velocity measurements, and actual orientations of some leaders.

\begin{Problem}[Orientation Localization]\label{prb:orientLocalization} 
Considering a \textit{twin-leader-follower} network of $n$ agents, under Assumptions \ref{ass:bidirection_sensing}-\ref{ass:non_colocated}, compute $\hat{\m{R}}_{i}$ for each follower $i,~i=3,\ldots,n,$ based on the directional measurements $(\m{b}_{ij}^i,\m{b}_{ji}^j)$, estimated orientations of its neighbors $\hat{\m{R}}_j,j \in \mc{N}_i,$ and the knowledge of the true orientations of the two or more leaders, i.e., $\m{R}_k \in SO(3),k \in \mc{V}_l$.
\end{Problem}

The second problem investigated in this work is to determine the locations of agents.
\begin{Problem}[Position Localization]\label{prb:orientEst}
Consider a \textit{twin-leader-follower} network of $n$ non-translating but possibly rotating agents with at least two leaders. Under Assumptions \ref{ass:bidirection_sensing}-\ref{ass:non_colocated}, for each follower $i$, determine its actual position, $\m{p}_i\in\mb{R}^{3}$, based on the estimate $\hat{\m{R}}_i$, the direction constraints $\m{b}_{ij}^i,j\in \mc{N}_i$, and absolute positions of some leaders, i.e., $\m{p}_k\in\mb{R}^{3},k \in \mc{V}_l$. 
\end{Problem}

\section{Orientation Localization}\label{sec:orient_local}
In this section, we present a differential equation constituting a continuous-time orientation localization law in $SO(3)$ that computes time-varying orientations of agents simultaneously using continuous-time directional vectors to multiple neighboring agents, angular velocity measurements, and actual orientations of some leaders. Further, the equilibrium set of the differential equation is first characterized and almost global asymptotically convergence of the estimated orientations is established.
\subsection{Error function and critical points}
Consider an agent $i\in \mc{V}_f$ which senses the local directions, $\m{b}_{ij}^i\in \mb{R}^3,$ to its neighboring agents $j\in \mc{N}_i$. If $|\mc{N}_i| =2$\footnote{When $|\mc{N}_i| >2$, agent $i$ and its neighbors are almost surely non-coplanar if they are placed randomly in $\mb{R}^3$.}, the third direction constraint is defined by the normalized cross product of the first two directions, for positive definiteness of $\m{K}_i$ in \eqref{eq:error_function_rewrite}.
The objective is to find an estimate, $\hat{\m{R}}_i\in SO(3)$, of the true orientation, $\m{R}_i$, that is the critical point of the following error function
\begin{align*}\label{eq:orient_cost_function}
\Phi_i(\hat{\m{R}}_i,\m{R}_i)  &={1}/{2}\textstyle\sum_{j\in \mc{N}_i}k_{ij}||\hat{\m{R}}_i\m{b}_{ij}^i-\m{b}_{ij}||^2\\
&=\textstyle\sum_{j\in \mc{N}_i}k_{ij}(1-\hat{\m{R}}_i\m{b}_{ij}^i\cdot\m{b}_{ij}), \numberthis
\end{align*}
which is sum of squared norms of all direction constraint errors. We do not assert that $\Phi_i$ can be evaluated from the measurements, but we shall show that it can be minimized from the measurements. In \eqref{eq:orient_cost_function}, positive constant gains, $k_{ij}\in \mb{R}$, are used to impose different weights on error terms in the error function. The above configuration error function is in the form of Wahba's cost function \cite{Wahba1965} (or an alternative formulation of the Procrustes problem \cite{Gower2004}) and used for attitude tracking control \cite{TLee2015tac} or attitude estimation of a rigid body \cite{Mahony2008tac,Izadi2016auto}. In the sequel, we follow techniques similar to those in \cite[Chap. 11]{Bullo2005spr},\cite{TLee2015tac} to design our orientation localization law.

Let 
$\Phi_{ij}\triangleq 1-\hat{\m{R}}_i\m{b}_{ij}^i\cdot\m{b}_{ij}=1-\text{tr}(\hat{\m{R}}_i\m{b}_{ij}^i\m{b}_{ij}^\top)=1-\text{tr}(\hat{\m{R}}_i\m{R}_i^\top\m{b}_{ij}\m{b}_{ij}^\top),$ where we use \eqref{eq:trace_property_1} and $\m{b}_{ij}^i=\m{R}_i^\top\m{b}_{ij}$. Let $\tilde{\m{Q}}_i\triangleq\hat{\m{R}}_i\m{R}_i^\top$ and hence $\Phi_{ij}=1-\text{tr}(\tilde{\m{Q}}_i\m{b}_{ij}\m{b}_{ij}^\top)$. Consider a vector in the tangent space of $SO(3)$ at the point $\hat{\m{R}}_i$ (resp. $\m{R}_i$) as $\delta\hat{\m{R}}_i =\hat{\m{R}}_i\eta_i^\wedge,\eta_i\in \mb{R}^3$, (resp. $\delta{\m{R}}_i ={\m{R}}_i\zeta_i^\wedge,\zeta_i\in \mb{R}^3$) \cite{Bullo2005spr}.
\begin{Lemma}\label{lm:left_der_error_function}
The derivative of the error function $\Phi(\hat{\m{R}}_i)$ with respect to $\hat{\m{R}}_i$ (resp. ${\m{R}}_i$) along the direction of $\hat{\m{R}}_i\eta_i^\wedge$ (resp. ${\m{R}}_i\zeta_i^\wedge$) is given by
\begin{align*}
&\m{D}_{\hat{\m{R}}_i}\Phi_i(\hat{\m{R}}_i,\m{R}_i)\cdot\hat{\m{R}}_i\eta_i^\wedge=\eta_i^\top\textstyle\sum_{j\in \mc{N}_i}\m{e}_{ij},\\
\big(\text{resp. }&\m{D}_{{\m{R}}_i}\Phi_i(\hat{\m{R}}_i,\m{R}_i)\cdot{\m{R}}_i\zeta_i^\wedge=-\zeta_i^\top\textstyle\sum_{j\in \mc{N}_i}\m{e}_{ij}\big),
\end{align*}
where $\m{e}_{ij}\triangleq k_{ij}(\hat{\m{R}}_i^\top\m{b}_{ij}\times \m{b}_{ij}^i)\in \mb{R}^3$, $j=1,\ldots,|\mc{N}_i|$.
\end{Lemma}
\begin{pf}
See Appendix \ref{app:left_der_error_proof}.
\end{pf}

We now study the critical points of $\Phi_i(\hat{\m{R}}_i)$. To proceed, we rewrite the error function as 
\begin{align*}\label{eq:error_function_rewrite}
\Phi_i&=\textstyle\sum_{j\in \mc{N}_i}k_{ij}-\textstyle\sum_{j\in \mc{N}_i}\text{tr}(  k_{ij}   \hat{\m{R}}_i\m{R}_i^\top\m{b}_{ij}\m{b}_{ij}^\top)\\
&=\textstyle\sum_{j\in \mc{N}_i}k_{ij}-\text{tr}(\tilde{\m{Q}}_i\m{K}_i) \numberthis
\end{align*}
 where $\m{K}_i\triangleq\sum_{j\in \mc{N}_i}k_{ij}\m{b}_{ij}(\m{b}_{ij})^\top\in \mb{R}^{3\times 3}$. For almost all positive scalars $k_{ij}$ the matrix $\m{K}_i$ in \eqref{eq:error_function_rewrite} has distinct eigenvalues\footnote{i.e. the discriminant of the cubic polynomial $\text{det}(\lambda\m{I}_3-\m{K}_i)$ is positive \cite{Irving2004}. It is noted that the discriminant has real coefficients and is a polynomial of the scalars $k_{ij}$. Futhermore, since all eigenvalues of $\m{K}_i$ are real the discriminant is nonnegative \cite{Irving2004} and the set of $k_{ij}$ such that the discriminant is zero is of measure zero \cite{Caron2005}. Consequently, $\m{K}_i$ has distinct eigenvalues for almost all positive scalars $k_{ij}$.}, and hence $\Phi_i(\hat{\m{R}}_i)$ has isolated critical points (See Lemma \ref{lm:critical_points}). Agent $i$ can design such scalars locally because $\m{K}_i$ is similar to the matrix $\sum_{j\in \mc{N}_i}k_{ij}\m{b}_{ij}^i(\m{b}_{ij}^i)^\top$. Since $\text{Range}\Big(k_{ij}\m{b}_{ij}(\m{b}_{ij})^\top\Big)=\text{span}\{\m{b}_{ij}\}$, 
it can be verified that $\m{K}_i$ is positive definite if and only if $\{\m{b}_{ij}\}_{j\in \mc{N}_i}$ are non-coplanar. Thus, $\m{K}_i$ can be decomposed as $\m{K}_i=\m{UG}\m{U}^\top$ where $\m{G}=\text{diag}\{\lambda_k(\m{K}_i)\},\lambda_k(\m{K}_i)>0,~k=1,2,3$, and $\m{U}\in O(3)$. Also note that $\text{tr}(\m{G})=\text{tr}(\m{K}_i)=\text{tr}(\sum_{j\in \mc{N}_i}k_{ij}\m{b}_{ij}\m{b}_{ij}^\top)=\sum_{j\in \mc{N}_i}k_{ij}\m{b}_{ij}^\top\m{b}_{ij}=\sum_{j\in \mc{N}_i}k_{ij}$. Consequently, one has
\begin{align*}
\Phi_i&=\text{tr}(\m{G})-\text{tr}(\tilde{\m{Q}}_i\m{UG}\m{U}^\top)=\text{tr}(\m{G})-\text{tr}(\m{G}\m{U}^\top\tilde{\m{Q}}_i\m{U})\\
&=\text{tr}(\m{G}(\m{I}_3-\m{U}^\top\tilde{\m{Q}}_i\m{U})),
\end{align*}
whose critical points are give as follows.
\begin{Lemma}\cite[Prop. 11.31]{Bullo2005spr}\label{lm:critical_points}
Let $\m{G}$ be a diagonal matrix with distinct positive entries and $\m{U}\in O(3)$. Then, $\Phi_i(\tilde{\m{Q}}_i)=tr(\m{G}(\m{I}_3-\m{U}^\top\tilde{\m{Q}}_i\m{U}))$ has four critical points given by
\begin{equation*}\tilde{\m{Q}}_i\in \{\m{I}_3,\m{U}\m{D}_1\m{U}^\top,\m{U}\m{D}_2\m{U}^\top,\m{U}\m{D}_3\m{U}^\top\},
\end{equation*}
where $\m{D}_i=2[\m{I}_3]_i[\m{I}_3]_i^\top-\m{I}_3$ and $[\m{I}_3]_i$ is the $i$-th column vector of $\m{I}_3$.
\end{Lemma}

Those critical points are clearly isolated in which $\tilde{\m{Q}}_i=\hat{\m{R}}_i\m{R}_i^\top=\m{I}_3$ is the desired point and $\text{tr}(\tilde{\m{Q}}_i)=-1$ for the three undesired points.

\subsection{Orientation localization law}
We now propose orientation localization law for each follower agent $i$ as
\begin{equation}\label{eq:orient_est_law}
\dot{\hat{\m{R}}}_i=\hat{\m{R}}_i\Omega_i^\wedge
\end{equation}
where the control vector $\Omega_i\in \mb{R}^3$ will be designed later and $\hat{\m{R}}_i(0)$ is initialized arbitrarily in $SO(3)$. Let $\tilde{\Omega}_i\triangleq\omega_i^i-\Omega_i$; we have the following lemma.
\begin{Lemma}\label{lm:time_der_error_function}
The vector, $\m{e}_i\triangleq\sum_{j\in \mc{N}_i}\m{e}_{ij}$, and error function, $\Phi_i$, in \eqref{eq:orient_cost_function} satisfy the following properties
\begin{enumerate}[(i)]
\item $||\dot{\m{e}}_i||\leq \sum_{j\in \mc{N}_i}k_{ij}||\tilde{\Omega}_i||+\bar{\omega}_i||\m{e}_i||$, where the positive constant $\bar{\omega}_i>0$ satisfies $||\omega_i||\leq\bar{\omega}_i$,
\item $\dot{\Phi}_i(\hat{\m{R}}_i, \m{R}_i)=-\tilde{\Omega}_i\cdot \m{e}_i,$
\item There exist constants $\sigma_i, \gamma_i>0$ such that $\sigma_i||\m{e}_i||^2\leq\Phi_i(\hat{\m{R}}_i, \m{R}_i)\le \gamma_i||\m{e}_i||^2$, where the upper bound holds when ${\Phi}_i<2\min \{\lambda_1+\lambda_2,\lambda_1+\lambda_3,\lambda_2+\lambda_3\},~(\lambda_k=\lambda(\m{K}_i),k=1,2,3)$.
\end{enumerate}
\end{Lemma}
\begin{pf} See Appendix \ref{app:time_der_error_func_proof}.
\end{pf}
The control vector $\Omega_i=\omega_i^i-\tilde{\Omega}_i$, where $\tilde{\Omega}_i\in \mb{R}^3$ is designed as
\begin{equation}\label{eq:omega_i}
\dot{\tilde{\Omega}}_i= -k_\omega\tilde{\Omega}_i+\textstyle\sum_{j\in \mc{N}_i}k_{ij}(\hat{\m{R}}_i^\top\hat{\m{R}}_j\m{b}_{ij}^j\times \m{b}_{ij}^i),
\end{equation}
where $k_\omega>0$ is a positive constant. The orientation localization law \eqref{eq:orient_est_law}-\eqref{eq:omega_i} is distributed in the sense that only locally measured directional vectors, i.e., $\m{b}_{ij}^i$, and information communicated from neighboring agents, i.e., the estimate of direction in the global coordinate frame, $\hat{\m{R}}_j\m{b}_{ij}^j$, are utilized. Since the right hand side of \eqref{eq:omega_i} is linear in $\tilde{\Omega}_i$ and the second term is bounded, $\tilde{\Omega}_i$ is uniformly continuous in $t$.
\subsection{Stability and convergence analysis}
We rewrite \eqref{eq:omega_i} as
\begin{align*}
\dot{\tilde{\Omega}}_i&= -k_\omega\tilde{\Omega}_i+\textstyle\sum_{j\in \mc{N}_i}k_{ij}(\hat{\m{R}}_i^\top\m{b}_{ij}\times \m{b}_{ij}^i\\
&\quad+\hat{\m{R}}_i^\top(\hat{\m{R}}_j-\m{R}_j)\m{b}_{ij}^j\times \m{b}_{ij}^i)\\
&=-k_\omega\tilde{\Omega}_i+\m{e}_i+\m{h}_i(\hat{\m{R}}_j,t) \numberthis \label{eq:orient_cascade_system}
\end{align*}
where $\m{h}_i(\hat{\m{R}}_j,t)=\sum_{j\in \mc{N}_i}k_{ij}(\hat{\m{R}}_i^\top(\hat{\m{R}}_j-\m{R}_j)\m{b}_{ij}^j\times \m{b}_{ij}^i)$. Due to the cascade structure of the leader-follower system we prove the almost global convergence of the estimated orientations using an induction argument.
\subsubsection{The first follower} For the first follower, i.e., agent $3$, we have $\m{h}_3=\m{0}$. Thus,
\begin{equation}\label{eq:orient_est_law_3}
\dot{\hat{\m{R}}}_3=\hat{\m{R}}_3(\omega_3^3-\tilde{\Omega}_3)^\wedge, ~\dot{\tilde{\Omega}}_3=-k_\omega\tilde{\Omega}_3+\m{e}_3.
\end{equation}
\begin{Theorem}
Suppose that Assumptions \ref{ass:bidirection_sensing}-\ref{ass:non_colocated} hold. Then, under the orientation localization law \eqref{eq:orient_est_law_3}, we have:
\begin{enumerate}[(i)]\label{thm:orient_agent_3}
\item The equilibrium points of \eqref{eq:orient_est_law_3} are given as $\{(\tilde{\m{Q}}_3,\tilde{\Omega}_3)~| ~\tilde{\m{Q}}_3\in \{\m{I}_3,\m{U}\m{D}_1\m{U}^\top,\m{U}\m{D}_2\m{U}^\top,\m{U}\m{D}_3\m{U}^\top\},$ $\tilde{\Omega}_3=\m{0}\}$, where $\m{D}_i$ and $\m{U}$ are defined in Lemma \ref{lm:critical_points}.
\item The desired equilibrium, $(\tilde{\m{Q}}_3=\m{I}_3,\tilde{\Omega}_3=\m{0})$ is almost globally asymptotically stable, $\tilde{\m{Q}}_3=\m{I}_3$ is the global minimum of $\Phi_3(\tilde{\m{Q}}_3)$ and the three undesired equilibria are unstable.
\end{enumerate}
\end{Theorem}
\begin{pf}
See Appendix \ref{app:orient_agent_3_proof}.
\end{pf}
It follows from the above theorem that $\hat{\m{R}}_3$ almost globally asymptotically converges to the true orientation $\m{R}_3$ as $t\rightarrow\infty$. For induction, we suppose that the corresponding result holds for agents $k-1,k-1\geq 3$, i.e., $\hat{\m{R}}_{k-1}\rightarrow\m{R}_{k-1}$ as $t\rightarrow\infty$ almost globally. We show that it is also true for the agent $k$ as follows.
\subsubsection{Follower $k$}
Using \eqref{eq:orient_est_law} and \eqref{eq:orient_cascade_system}, we have
\begin{equation}\label{eq:orient_est_law_k}
\dot{\hat{\m{R}}}_{k}=\hat{\m{R}}_{k}(\omega_k^k-\tilde{\Omega}_{k})^\wedge, ~\dot{\tilde{\Omega}}_{k}=-k_\omega\tilde{\Omega}_{k}+\m{e}_{k}+\m{h}_{k}(t),
\end{equation}
where $\m{h}_{k}(\hat{\m{R}}_{k},t)=\sum_{j\in \mc{N}_{k}}k_{kj}(\hat{\m{R}}_k^\top(\hat{\m{R}}_j-\m{R}_j)\m{b}_{kj}^j\times \m{b}_{kj}^k)$ which is clearly bounded and converges to zero asymptotically since $\hat{\m{R}}_{j}\rightarrow\m{R}_{j},~\forall j=1,\ldots,k-1$. Note that $\m{h}_{k}(\hat{\m{R}}_{k},t)$ can be considered as an additive input to the nominal system
\begin{equation}\label{eq:orient_nominal_system_k}
\dot{\hat{\m{R}}}_{k}=\hat{\m{R}}_{k}(\omega_k^k-\tilde{\Omega}_{k})^\wedge, ~\dot{\tilde{\Omega}}_{k}=-k_\omega\tilde{\Omega}_{k}+\m{e}_{k}.
\end{equation}
It is noted that the above system is in a similar form to \eqref{eq:orient_est_law_3} and hence the following result follows directly.
\begin{Lemma}\label{lm:orient_agent_k}
Consider the nominal system \eqref{eq:orient_nominal_system_k} under the Assumptions \ref{ass:bidirection_sensing}-\ref{ass:non_colocated}, then:
\begin{enumerate}[(i)]
\item The four isolated equilibrium points of \eqref{eq:orient_nominal_system_k} are given as $\{(\tilde{\m{Q}}_k,\tilde{\Omega}_k)~| ~\tilde{\m{Q}}_k\in \{\m{I}_3,\m{U}\m{D}_1\m{U}^\top,\m{U}\m{D}_2\m{U}^\top,$ $\m{U}\m{D}_3\m{U}^\top\},\tilde{\Omega}_{k}=\m{0}\}$, where $\m{D}_i$ and $\m{U}$ are defined in Lemma \ref{lm:critical_points}.
\item The desired equilibrium, $(\tilde{\m{Q}}_k=\m{I}_3,\tilde{\Omega}_k=\m{0})$ is almost globally asymptotically stable while the three undesired equilibria are unstable.
\end{enumerate}
\end{Lemma}
The perturbed system \eqref{eq:orient_est_law_k} is linear in $\tilde{\Omega}_k$ and $\m{e}_k+\m{h}_k$ is bounded. Thus $\tilde{\Omega}_k$ is bounded and the perturbed system \eqref{eq:orient_est_law_k} also satisfies the ultimate boundedness property and is input-to-state stable (ISS) with respect to the input $\m{h}_k(t)$ as will be shown below. Define the set $\mc{S}_k\triangleq\{\tilde{\m{Q}}_k~|~\Phi_k(\tilde{\m{Q}}_k)<\phi_k\}$, where $\phi_k=2\min \{\lambda_1+\lambda_2,\lambda_1+\lambda_3,\lambda_2+\lambda_3\},~\{\lambda_i\}_{i=1,2,3}=\lambda(\m{K}_k)$, or, i.e., the minimum value of $\Phi_k$ evaluated at the three undesired critical points. 
\begin{Lemma}\label{lm:ultimate_boundedness}
Suppose that Assumptions \ref{ass:bidirection_sensing}-\ref{ass:non_colocated} hold. The perturbed system \eqref{eq:orient_est_law_k} is almost input-to-state stable (ISS) with respect to $\m{h}_k(t)$. 
\end{Lemma}
\begin{pf}
See Appendix \ref{app:ultimate_boundedness}.
\end{pf}
\begin{Theorem}\label{thm:orient_agent_k_ISS}
Suppose that Assumptions \ref{ass:bidirection_sensing}-\ref{ass:non_colocated} hold. Then, the desired equilibrium point, $(\hat{\m{R}}_k=\m{R}_k,\tilde{\Omega}_k=\m{0})$, of the system \eqref{eq:orient_est_law_k} is almost globally asymptotically stable. 
\end{Theorem}
\begin{pf}
First, the desired equilibrium point $(\hat{\m{R}}_k=\m{R}_k,\tilde{\Omega}_k=\m{0})$ of the unforced system \eqref{eq:orient_nominal_system_k} is almost globally asymptotically stable (Lemma \ref{lm:orient_agent_k}). The other undesired equilibria are isolated and unstable. The perturbed system \eqref{eq:orient_est_law_k} satisfies ultimate boundedness and is input-to-state stable w.r.t $\m{h}_{k}$ (Lemma \ref{lm:ultimate_boundedness}). The input $\m{h}_{k}(t)$ is bounded and vanishes asymptotically as $t\rightarrow\infty$.
Consequently, the desired equilibrium point, $(\hat{\m{R}}_k=\m{R}_k,\tilde{\Omega}_k=\m{0})$, of the system \eqref{eq:orient_est_law_k} is almost globally stable \cite{Angeli2011tac,Angeli2004tac}. \eprf
\end{pf} 
It follows that $\hat{\m{R}}_k(t)\rightarrow \m{R}_k$ almost globally asymptotically as $t\rightarrow \infty$. Finally, by invoking mathematical induction, the above theorem holds for all $k=3,\ldots,n$.
\begin{Corollary}\label{coroll:robust_to_noise}
Suppose direction measurements include bounded additive measurement noise. Then for a sufficiently small bound, $(\tilde{\m{Q}}_k,\tilde{\Omega}_k)$ converges to a neighborhood of the desired equilibrium $(\m{I}_3,\m{0})$ of \eqref{eq:orient_est_law_k}.
\end{Corollary}
\begin{pf}
See Appendix \ref{app:robust_to_noise}.
\end{pf}

\section{Position Localization} \label{sec:pos_local}
This section investigates the position localization problem. The aim of the position localization is to determine the positions of all followers using locally measured directions $\m{b}_{ij}^i$, the estimated orientation $\hat{\m{R}}_i$ of agent $i$ and the absolute positions of some leaders. For this, we first study the uniqueness of the target positions of the followers and present a distributed localization law for each agent. Under the proposed position localization law, estimated positions of all followers converge globally and asymptotically to the true positions.
\subsection{Unique target configuration}
The uniqueness of the target configuration (the actual positions of agents) is a key property of the network that allows us to localize the network. In the sequel, under the noncollocation and non-collinearity  conditions of the true positions of the agents in Assumption \ref{ass:non_colocated}, we show that the target configuration is uniquely defined using the direction constraints, estimate of orientation of agent $i$, and the absolute positions of some leaders. The following result is similar to \cite[Lm. 1]{Minh2017tac}.
\begin{Lemma}[Unique Target Configuration]\label{lm:unique_position}
Consider the twin-leader-follower network with two or more leaders and locally measured directions $\{\m{b}_{ij}^i\}_{(i,j)\in \mc{E}}$. Suppose that Assumptions \ref{ass:bidirection_sensing}-\ref{ass:non_colocated} hold, and the orientation of agent $i$, $\m{R}_i\in SO(3)$, is available to $i$ or otherwise can be estimated, e.g. Problem \ref{prb:orientEst}. Then the actual position of the agent $i$, ($i\geq3$), i.e., $\m{p}_i\in \mb{R}^3$ is uniquely determined from its direction constraints $\{\m{b}_{ij}^i\}_{j\in \mc{N}_i}$ and the positions of its neighbors $\{\m{p}_j\}_{j\in \mc{N}_i}$. Furthermore, $\m{p}_i$ is uniquely computed as
\begin{equation}\label{eq:unique_position}
\m{p}_i = \Big( \textstyle\sum_{j\in \mc{N}_i}\m{P}_{\m{b}_{ij}}\Big)^{-1}\sum_{j\in \mc{N}_i}\m{P}_{\m{b}_{ij}}\m{p}_j,
\end{equation}
where $\m{b}_{ij}=\m{R}_i\m{b}_{ij}^i$, and $\m{P}_{\m{b}_{ij}}\in \mb{R}^{3\times 3}$ denotes the projection matrix as defined in \eqref{eq:projection_matrix}. 
\end{Lemma}
\begin{pf}
The position $\m{p}_i$ of agent $i$ must satisfy the direction constraints
\begin{align*}
\m{P}_{\m{b}_{ij}}(\m{p}_i-\m{p}_j)&=\m{0},~\forall j\in \mc{N}_i.
\end{align*}
It follows from the above constraints that
\begin{equation}\label{eq:position_i_constraint}
\Big( \textstyle\sum_{j\in \mc{N}_i}\m{P}_{\m{b}_{ij}}\Big)\m{p}_i=\textstyle\sum_{j\in \mc{N}_i}\m{P}_{\m{b}_{ij}}\m{p}_j,
\end{equation}
Since $\text{null}(\m{P}_{\m{b}_{ij}})=\text{span}\{\m{b}_{ij}\}$ and position of follower $i$ is not collinear with two or more of its neighboring agents (Assumption \ref{ass:non_colocated}), we have ${\cap} \{\text{null}(\m{P}_{\m{b}_{ik}})\}_{j\in \mc{N}_i}=\m{0}$.
As a result, $ \sum_{j\in \mc{N}_i}\m{P}_{\m{b}_{ij}}$ is positive definite and hence invertible. Thus, $\m{p}_i$ is uniquely computed as \eqref{eq:unique_position}. \eprf
\end{pf}

\begin{Remark}
It is worth noting that $\m{R}_i\in SO(3)$ and $\{\m{p}_j\}_{j\in \mc{N}_i}$ are not available to $i$ initially but these quantities do become available, i.e., when the corresponding quantities for its neighbors have been computed and made available to it. In the following subsection, we present a position localization law which runs in parallel with the aforementioned orientation estimation scheme \eqref{eq:orient_est_law}.
\end{Remark}

\subsection{Proposed position localization law}
Each follower agent $i$ holds an initial estimate of its position $\hat{\m{p}}_i(0)\in \mb{R}^3$. For each follower $i$, we propose the following position localization law
\begin{equation}\label{eq:pos_local_law}
\dot{\hat{\m{p}}}_i=-\hat{\m{R}}_i\textstyle\sum_{j\in \mc{N}_i}k_{p_{ij}}\m{P}_{\m{b}_{ij}^i}\hat{\m{R}}_i^\top(\hat{\m{p}}_i-\hat{\m{p}}_j),
\end{equation}
where, $k_{p_{ij}}>0$ is a positive gain, the local projection matrix $\m{P}_{\m{b}_{ij}^i}=\m{I}_3-\m{b}_{ij}^i(\m{b}_{ij}^i)^\top=\m{R}_i^\top(\m{I}_3-\m{b}_{ij}\m{b}_{ij}^\top)\m{R}_i=\m{R}_i^\top\m{P}_{\m{b}_{ij}}\m{R}_i$, and $\hat{\m{p}}_i(0)$ is initialized arbitrarily.
The localization law \eqref{eq:pos_local_law} is implemented using only local direction measurements $\m{b}_{ij}^i$, estimate of orientation $\hat{\m{R}}_i$, and estimates of its neighbors' positions $\hat{\m{p}}_j$ which are communicated by agents $j\in\mc{N}_i$ (in the case of leaders, $\hat{\m{p}}_i=\m{p}_i,\forall i \in \mc{V}_l$). The estimation law \eqref{eq:pos_local_law} is linear in the estimated state $\hat{\m{p}}(t)$, thus, and so the right side is globally Lipschitz in $\hat{\m{p}}(t)$.
\begin{Remark}
Given absolute positions of some leaders $\m{p}_i\in \mb{R}^3,i\in \mc{V}_l$, where $\mc{V}_l$ denotes the leader set, and since (in \eqref{eq:pos_local_law}) $\hat{\m{R}}_i\rightarrow\m{R}_i$ as $t\rightarrow\infty$, the steady-state solutions (if they exist) to \eqref{eq:pos_local_law} satisfy the following direction constraints
\begin{equation}\label{eq:position_constraints}
\begin{cases} \m{P}_{\m{b}_{ij}}(\hat{\m{p}}_i-\hat{\m{p}}_j)=0, &\mbox{} \forall (i,j)\in \mc{E},\\ 
\hat{\m{p}}_i=\m{p}_i & \mbox{} \forall i\in \mc{V}_l. 
\end{cases}
\end{equation}
From Lemma \eqref{lm:unique_position} and due to the cascade structure of the system, it can be shown that the true network location, $\m{p}\in \mb{R}^{3n}$, is the unique solution to \eqref{eq:position_constraints}. 
\end{Remark}
\subsection{Analysis}
We rewrite the localization law \eqref{eq:pos_local_law} as follows
\begin{align*}
\dot{\hat{\m{p}}}_i=\m{f}_i(\hat{\m{p}},t) -\m{h}_i(\hat{\m{p}},\hat{\m{R}})
\end{align*}
where 
$\m{f}_i(\hat{\m{p}}):=-{\m{R}}_i\textstyle\sum_{j\in \mc{N}_i}k_{p_{ij}}\m{P}_{\m{b}_{ij}^i}{\m{R}}_i^\top(\hat{\m{p}}_i-{\m{p}}_j)\\
=-\textstyle\sum_{j\in \mc{N}_i}k_{p_{ij}}\m{P}_{\m{b}_{ij}}(\hat{\m{p}}_i-{\m{p}}_j)$ and
\begin{align*}\m{h}_i(\hat{\m{p}},\hat{\m{R}})&:=-(\hat{\m{R}}_i-{\m{R}}_i)\textstyle\sum_{j\in \mc{N}_i}k_{p_{ij}}\m{P}_{\m{b}_{ij}^i}{\m{R}}_i^\top(\hat{\m{p}}_i-\hat{\m{p}}_j)\\
&-\hat{\m{R}}_i\textstyle\sum_{j\in \mc{N}_i}k_{p_{ij}}\m{P}_{\m{b}_{ij}^i}(\hat{\m{R}}_i^\top-{\m{R}}_i^\top)(\hat{\m{p}}_i-\hat{\m{p}}_j)\\
&-{\m{R}}_i\textstyle\sum_{j\in \mc{N}_i}k_{p_{ij}}\m{P}_{\m{b}_{ij}^i}{\m{R}}_i^\top({\m{p}}_j-\hat{\m{p}}_j).
\end{align*}
The above dynamics can be written in a more compact form
\begin{equation}\label{eq:pos_local_compact_form}
\dot{\hat{\m{p}}}=\m{f}(\hat{\m{p}},t)+\m{h}(\hat{\m{p}},\hat{\m{R}}),
\end{equation}
where the stack vectors $\m{f}(\hat{\m{p}})=[\m{f}_1,\m{f}_2,\ldots,\m{f}_n]^\top$ and $\m{h}(\hat{\m{p}},t)=[\m{h}_1,\m{h}_2,\ldots,\m{h}_n]^\top$. Due to the cascade structure of the system \eqref{eq:pos_local_compact_form}, we will study \eqref{eq:pos_local_compact_form} using the stability theory for cascade systems \cite{Angeli2004tac}.
\subsection{Stability Analysis}
Consider $\m{h}(\hat{\m{p}},t)$ in \eqref{eq:pos_local_compact_form} as an input to the following nominal system
\begin{equation}\label{eq:pos_nominal_system}
\dot{\hat{\m{p}}}=\m{f}(\hat{\m{p}}).
\end{equation}
The boundedness of the estimates of positions is provided in the following lemma.
\begin{Lemma}
Under Assumptions \ref{ass:bidirection_sensing}-\ref{ass:non_colocated}, the cascade system \eqref{eq:pos_local_compact_form} satisfies the ultimate boundedness property. That is, the estimates $\hat{\m{p}}_i ~(i=3,\ldots,n)$  are bounded for all time $t>0$.

\end{Lemma}
\begin{Lemma}[Global Asymptotic Stability]\label{lm:pos_nominal_system}
Under Assumptions \ref{ass:bidirection_sensing}-\ref{ass:non_colocated}, the desired equilibrium ${\hat{\m{p}}}=\m{p}$ of the nominal system \eqref{eq:pos_nominal_system} is globally exponentially stable. 
\end{Lemma}
\begin{pf}
For each follower $i\in \mc{V}_f$, the equilibrium of $\dot{\hat{\m{p}}}_i=\m{f}(\hat{\m{p}}_i)$ satisfies $\m{f}(\hat{\m{p}}_i)=\m{0}\Leftrightarrow \sum_{j\in \mc{N}_i}k_{p_{ij}}\m{P}_{\m{b}_{ij}}(\hat{\m{p}}_i-{\m{p}}_j)=\m{0}\Leftrightarrow \Big( \sum_{j\in \mc{N}_i}k_{p_{ij}}\m{P}_{\m{b}_{ij}}\Big)\hat{\m{p}}_i=\sum_{j\in \mc{N}_i}k_{p_{ij}}\m{P}_{\m{b}_{ij}}\m{p}_j$. Since agent $i$ is not colinear with two or more of its neighbors (Assumption \ref{ass:non_colocated}), $\hat{\m{p}}_i=\m{p}_i$ is the unique solution to the equation (Lemma \ref{lm:unique_position}). Consequently, ${\hat{\m{p}}}=\m{p}$, is the unique equilibrium of the nominal system \eqref{eq:pos_nominal_system}.

Consider a Lyapunov function $V_i=1/2(\hat{\m{p}}_i-\m{p}_i)^2$, which is positive definite, continuously differentiable, and radially unbounded. The derivative of $V_i$ along the trajectory of \eqref{eq:pos_nominal_system} is given as
\begin{align*}
\dot{V}_i(t)&= (\hat{\m{p}}_i-\m{p}_i)^\top\dot{\hat{\m{p}}}_i\\
&=-(\hat{\m{p}}_i-\m{p}_i)^\top\sum_{j\in \mc{N}_i}k_{p_{ij}}\m{P}_{\m{b}_{ij}}(\hat{\m{p}}_i-{\m{p}}_j)\\
&=-(\hat{\m{p}}_i-\m{p}_i)^\top\sum_{j\in \mc{N}_i}k_{p_{ij}}\m{P}_{\m{b}_{ij}}(\hat{\m{p}}_i-{\m{p}}_i+{\m{p}}_i-{\m{p}}_j)\\
&=-(\hat{\m{p}}_i-\m{p}_i)^\top\big(\sum_{j\in \mc{N}_i}k_{p_{ij}}\m{P}_{\m{b}_{ij}}\big)(\hat{\m{p}}_i-{\m{p}}_i),
\end{align*}
where the last equality follows from $\m{P}_{\m{b}_{ij}}({\m{p}}_i-{\m{p}}_j)=\m{0}$, for all $j\in \mc{N}_i$. Since $\big(\sum_{j\in \mc{N}_i}k_{p_{ij}}\m{P}_{\m{b}_{ij}}\big)$ is positive definite (Lemma \ref{lm:unique_position}), $\dot{V}_i(t)$ is negative definite. This completes the proof.\eprf
\end{pf}
\begin{Theorem}[Input-to-state Stability]
Under Assumptions \ref{ass:bidirection_sensing}-\ref{ass:non_colocated}, the cascade system \eqref{eq:pos_local_compact_form} is input-to-state stable with respect to the input $\m{h}(\hat{\m{p}},\hat{\m{R}})$. Further, ${\hat{\m{p}}}(t)\rightarrow\m{p}$ almost globally and asymptotically as $t\rightarrow\infty$. 
\end{Theorem}
\begin{pf}
Due to the cascade structure of the twin-leader-follower network we provide a proof by using mathematical induction. The almost global asymptotic convergence of the localized position of the first follower follows directly since the desired equilibrium, $\hat{\m{p}}_3=\m{p}_3$, of the nominal system, $\dot{\hat{\m{p}}}_3=\m{f}_3(\hat{\m{p}}_3,t)$, is globally exponentially stable (Lemma \ref{lm:pos_nominal_system}), and the input is bounded and $\m{h}_3(t)\rightarrow 0$ asymptotically due to $\hat{\m{R}}_i\rightarrow\m{R}_i$ almost globally as $t\rightarrow \infty$ and $\hat{\m{p}}_j=\m{p}_j,\forall j\in \mc{N}_3$. It can be shown similarly for all other followers using the facts that the convergence of estimated position of an follower is not influenced by the latter agents in the network and the orientations and positions of earlier agents converge to the actual poses asymptotically. This completes the proof.\eprf
\end{pf}
\begin{Remark}
When the two leaders $1$ and $2$ do not have knowledge of their actual poses (and there are no other redundant leaders), leader $1$ can fix the translation and rotation ambiguities in the pose estimation by choosing an arbitrary guess of its pose, e.g., $(\hat{\m{R}}_1,\hat{\m{p}}_1)\in SO(3)\times \mb{R}^3$, while leader $2$ can determine the scale factor by selecting an arbitrary guess of the distance to agent $1$, e.g., $\hat{d}_{21}>0$ (see e.g., \cite{Minh2017tac}). In particular, the leader $2$ can easily compute relative orientation to leader $1$, $\m{R}_{21}\in SO(3)$ using the directions between two leaders and the direction of each leader to the first follower $3$ (or i.e., based the triangulation network of the agents $(1,2,3)$, see e.g., \cite{Quoc2018cdc}). Then, the leader $2$ can use $(\hat{\m{R}}_2=\hat{\m{R}}_1\m{R}_{21}^\top,\hat{\m{p}}_2=\hat{\m{p}}_1-\hat{d}_{21}\hat{\m{R}}_2\m{b}_{21}^{2})$ as the guess of its pose, and consequently, the poses of the followers are determined up to a translation $\hat{\m{p}}_1$, a rotation $\hat{\m{R}}_1$, and a scale $\hat{d}_{21}/d_{21}$, where $d_{21}$ is the actual distance between two leaders.
\end{Remark}
\section{Simulation}\label{sec:Sim}
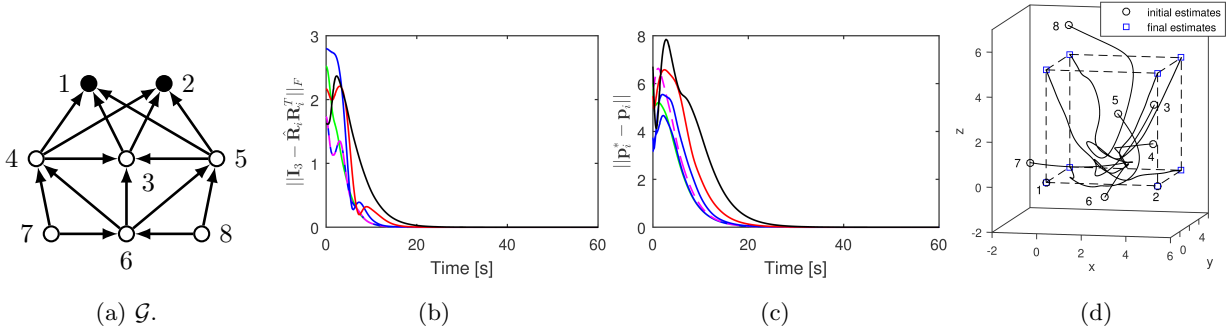
\begin{figure*}[t]
\centering
\begin{subfigure}[t]{0.2\textwidth}
\centering
\begin{tikzpicture}[scale=1.]
\node[place,black] (p1) at (0,0.5) [label=left:$1$]{};
\node[place,black] (p2) at (1,0.5) [label=right:$2$]{};
\node[place] (p3) at (0.5,-0.5) [label=below right:$3$]{};
\node[place] (p4) at (-0.7,-0.5) [label=left:$4$]{};
\node[place] (p5) at (1.7,-0.5) [label=right:$5$]{};
\node[place] (p6) at (.5,-1.5) [ label=below:$6$]{};
\node[place] (p7) at (-.5,-1.5) [label=left:$7$]{};
\node[place] (p8) at (1.5,-1.5) [label=right:$8$]{};
\draw[line width=1pt,->] (p3)[frame]  -- (p1);
\draw[line width=1pt,->] (p3)[frame]  -- (p2);
\draw[line width=1pt,->] (p4)[frame]  -- (p3);
\draw[line width=1pt,->] (p4)[frame]  -- (p1);
\draw[line width=1pt,->] (p5)[frame]  -- (p1);
\draw[line width=1pt,->] (p5)[frame]  -- (p3);
\draw[line width=1pt,->] (p5)[frame]  -- (p2);
\draw[line width=1pt,->] (p4)[frame]  -- (p2);
\draw[line width=1pt,->] (p6)[frame]  -- (p4);
\draw[line width=1pt,->] (p6)[frame]  -- (p5);
\draw[line width=1pt,->] (p6)[frame]  -- (p3);
\draw[line width=1pt,->] (p8)[frame]  -- (p6);
\draw[line width=1pt,->] (p8)[frame]  -- (p5);
\draw[line width=1pt,->] (p7)[frame]  -- (p4);
\draw[line width=1pt,->] (p7)[frame]  -- (p6);
\end{tikzpicture}
\caption{$\mc{G}$.}
\label{fig:sim_graph}
\end{subfigure}
\begin{subfigure}[t]{0.25\textwidth}
\includegraphics[height=3.3cm]{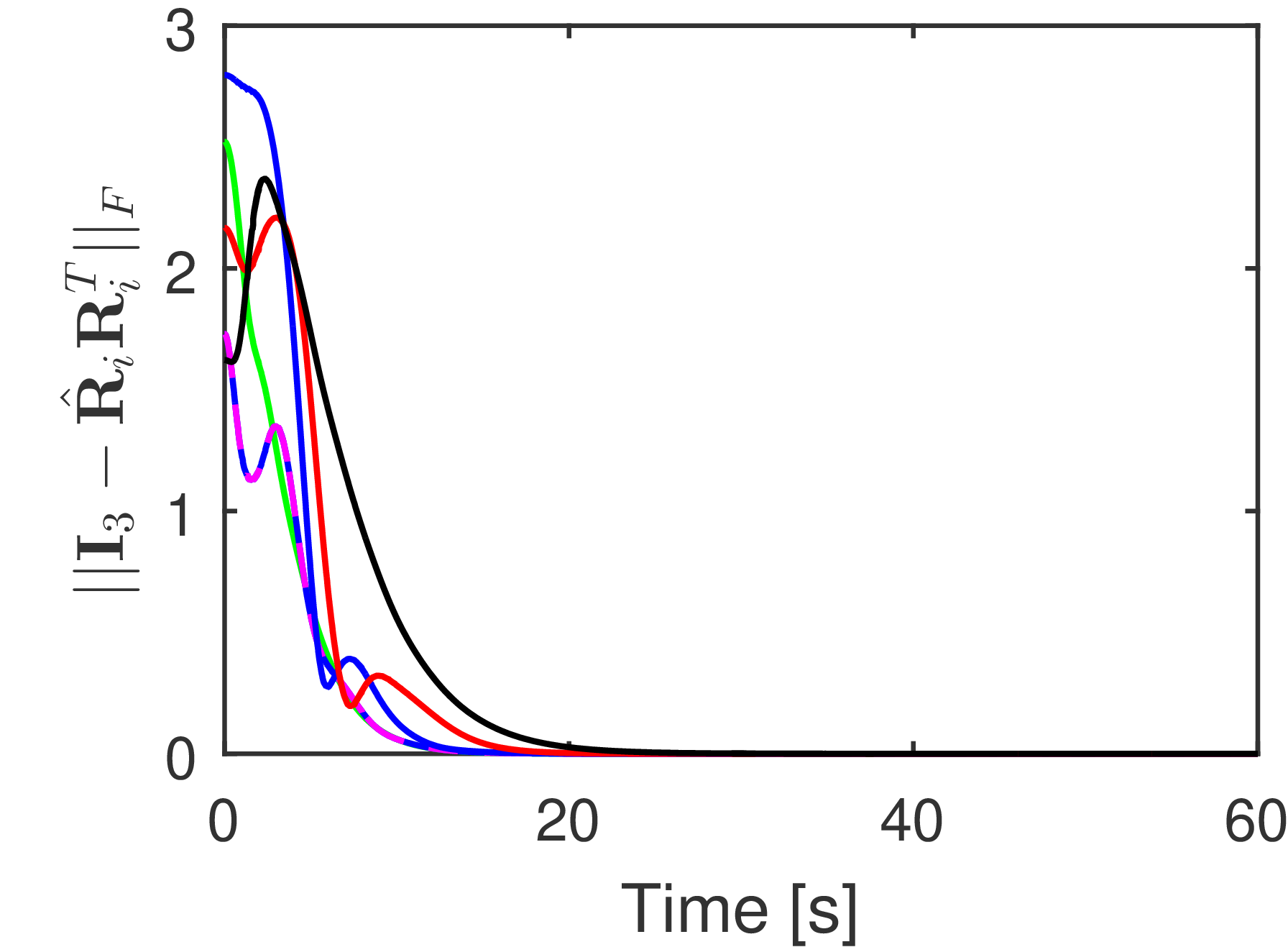}
\caption{}
\label{fig:sim_orient_err}
\end{subfigure}
\begin{subfigure}[t]{0.25\textwidth}
\includegraphics[height=3.3cm]{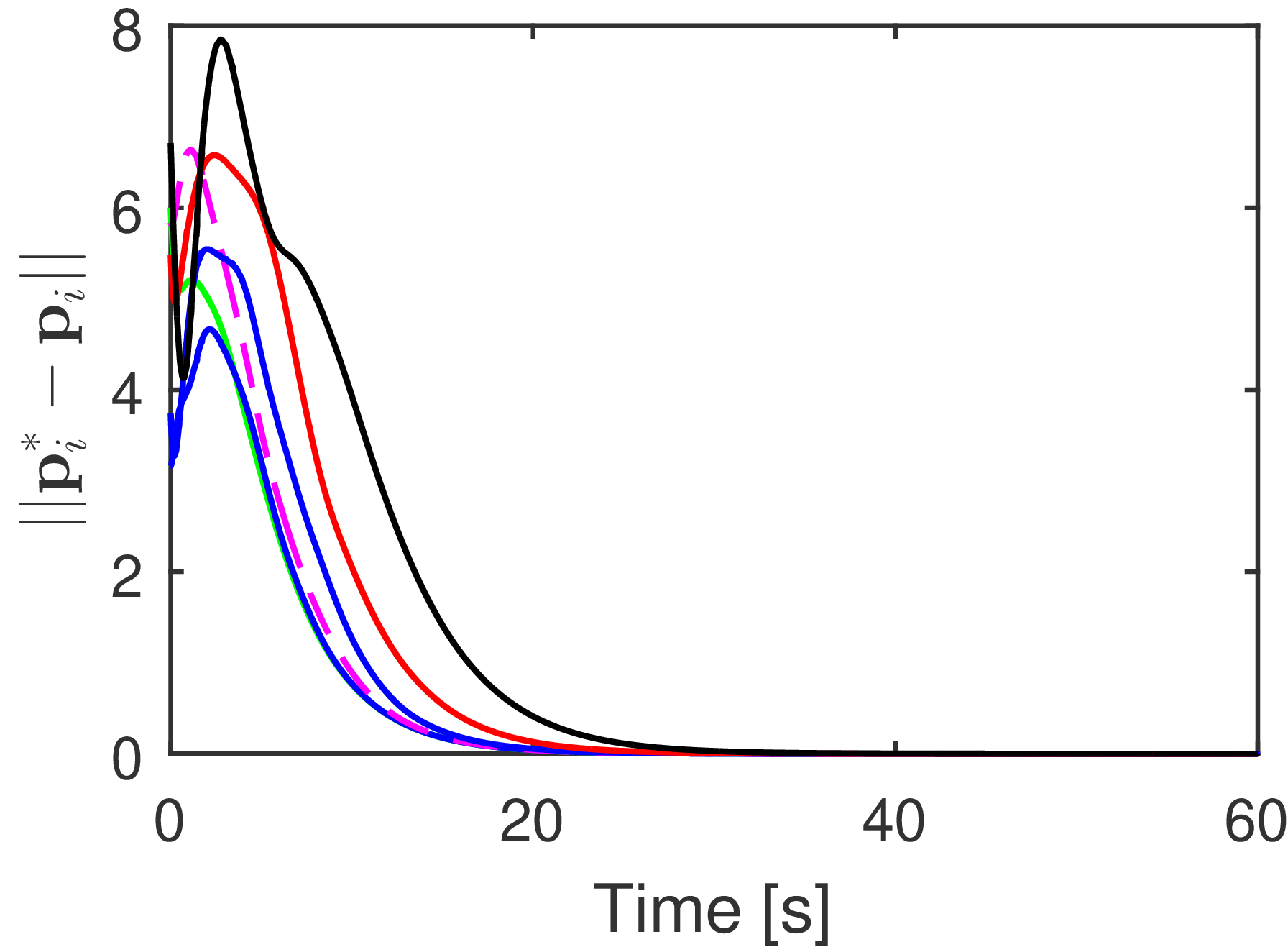}
\caption{}
\label{fig:sim_pos_err}
\end{subfigure}
\begin{subfigure}[t]{0.22\textwidth}
\includegraphics[height=4cm]{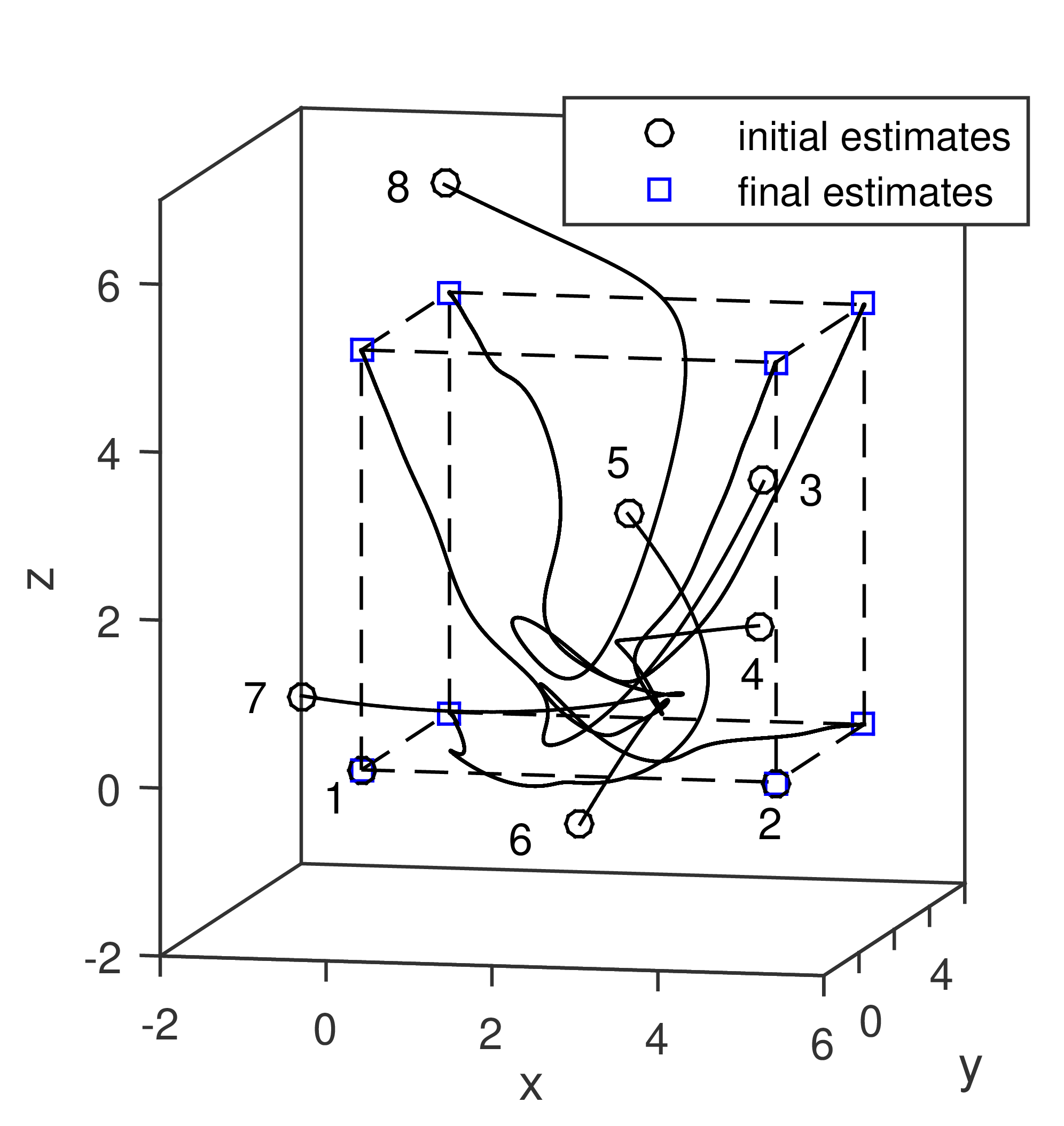}
\caption{}
\label{fig:sim_trajectory}
\end{subfigure}
\caption{Pose localization of six followers in $\mb{R}^3$ with accurate direction measurements. (a) A \textit{twin-leader-follower} network for simulation (b) Orientation induced norm errors. (c) Position localization norm errors. (d) Evolutions of estimated positions of all agents.}
\label{fig:sim}
\end{figure*}

Consider a twin-leader-follower network of eight agents in $\mb{R}^3$ whose graph topology is given in Fig. \ref{fig:sim_graph}. Agents $1$ and $2$ are leaders which know their actual poses. The initial and final estimates of network location are shown in Fig. \ref{fig:sim_trajectory}. The first followers $3,4,$ and $5$ keep rotating about their local $x,y$, and $z$ axes with the same angular velocity of $0.15$ rad/s, i.e., $\omega_3^3=[0.15,0,0]^\top,\omega_4^4=[0,0.15,0]^\top$ and $\omega_5^5=[0,0,0.15]^\top$, respectively; orientations of the other agents are fixed. The initial orientations of agents are chosen randomly. The actual configuration of the network is a cube of side length of $5$.

Simulation results are provided in Fig. \ref{fig:sim}. It can be seen that the estimated poses of the agents converge to the actual poses asymptotically as  the orientation induced norm errors, $||\m{I}_3-\hat{\m{R}}_i\m{R}_i^T||_F$, where $||\m{A}||_F=\sqrt{\text{tr}(\m{A}^T\m{A})}$ denotes the Frobenius norm of a matrix $\m{A}$, and position norm errors, $||\hat{\m{p}}_i-\m{p}_i||$, converge to zero asymptotically (See Figs. \ref{fig:sim_orient_err} and \ref{fig:sim_pos_err}, respectively).
\subsection{Pose Localization under direction measurement noise}
\begin{figure*}[t]
\centering
\begin{subfigure}[t]{0.3\textwidth}
\includegraphics[height=3.3cm]{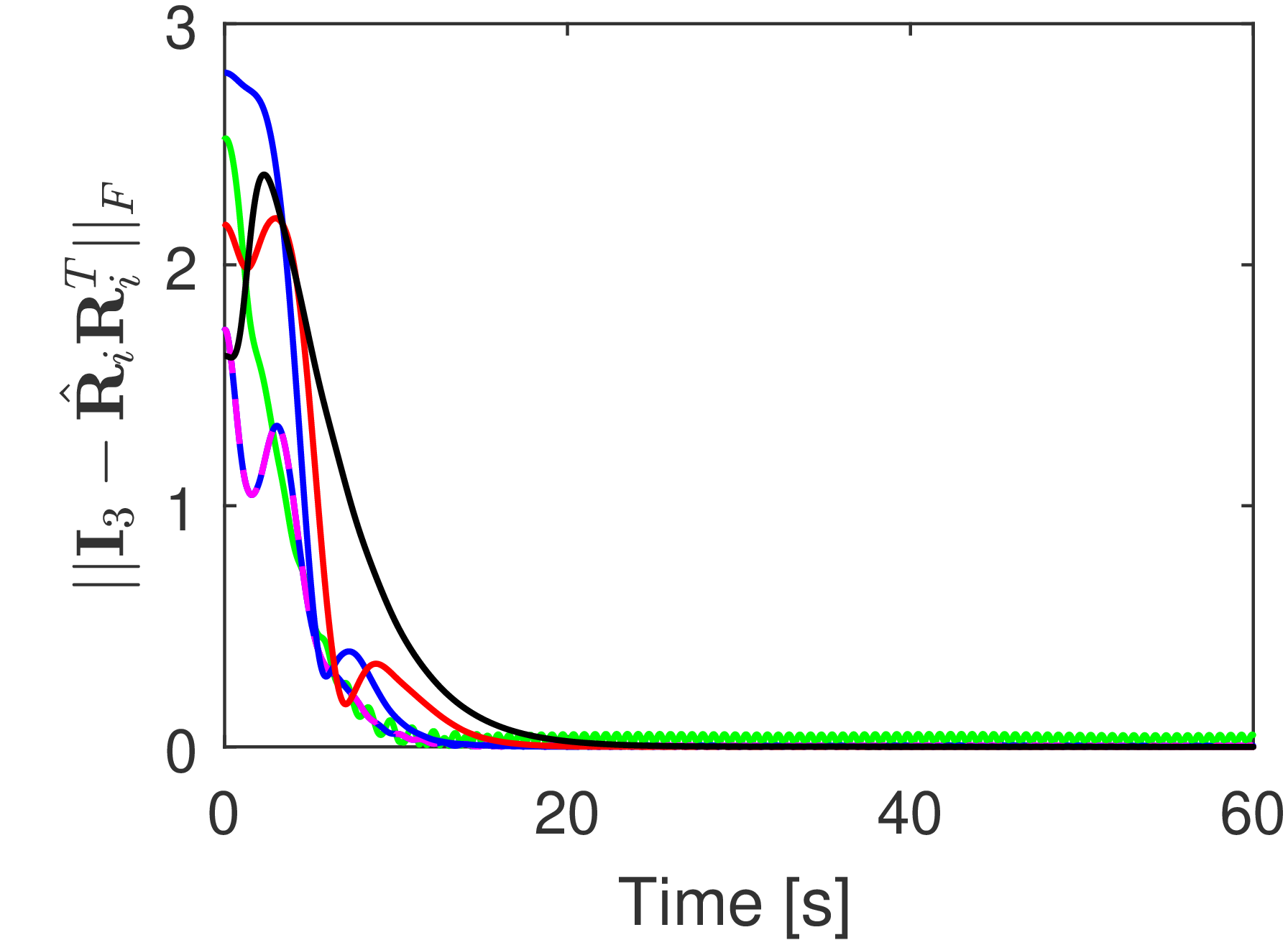}
\caption{}
\label{fig:sim_orient_err_noise}
\end{subfigure}
\begin{subfigure}[t]{0.3\textwidth}
\includegraphics[height=3.3cm]{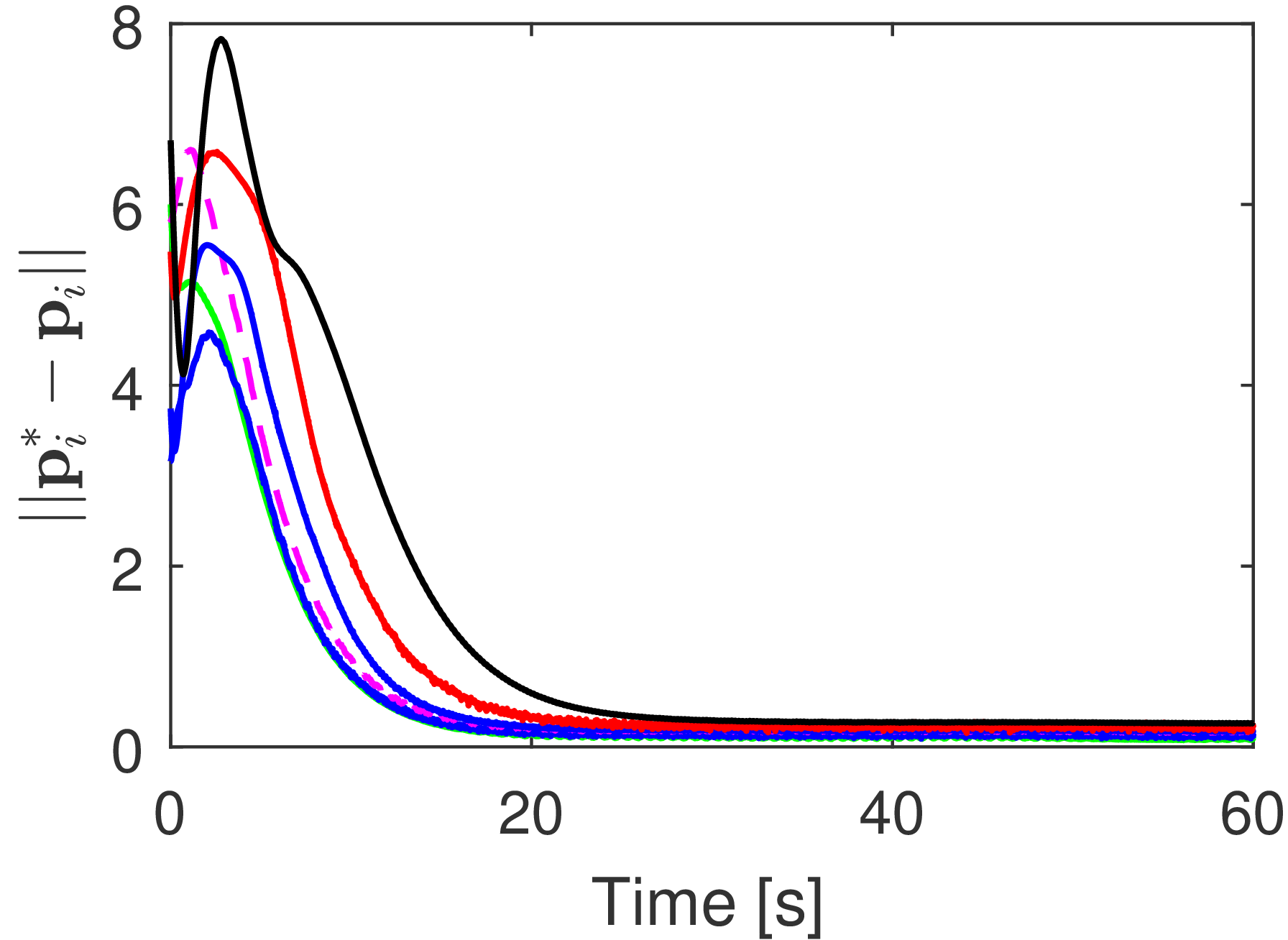}
\caption{}
\label{fig:sim_pos_err_noise}
\end{subfigure}
\begin{subfigure}[t]{0.3\textwidth}
\includegraphics[height=4cm]{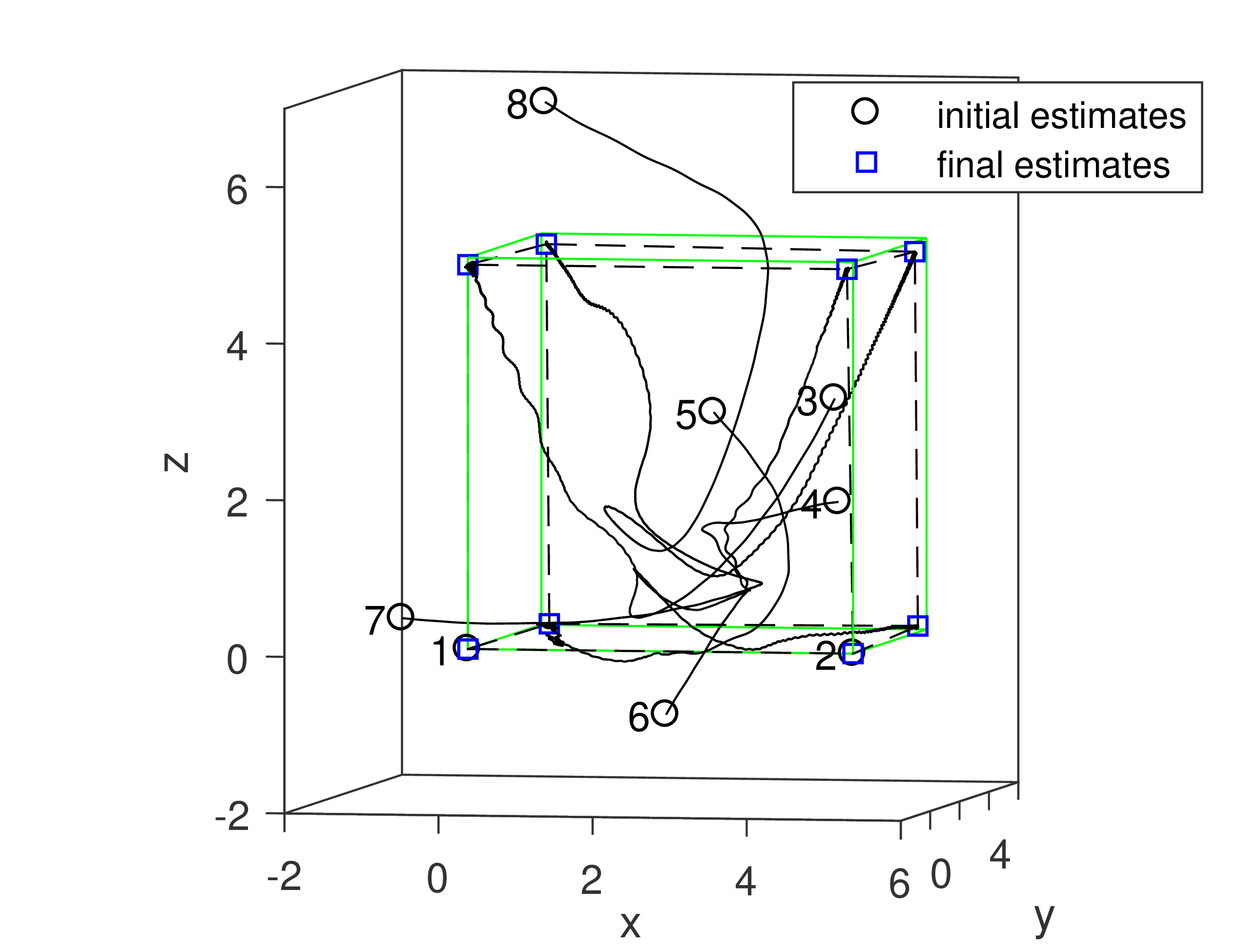}
\caption{}
\label{fig:sim_trajectory_noise}
\end{subfigure}
\caption{Pose localization of six followers in $\mb{R}^3$ under direction measurement noise. (a) Orientation induced norm errors. (b) Position localization norm errors. (c) Evolutions of estimated positions of all agents and the desired positions (green box).}
\label{fig:sim_noise}
\end{figure*}

We now assume that each true direction measurement $\m{b}_{ij}\in \mb{R}^3$ is contaminated by noise as follows
\begin{equation*}
\tilde{\m{b}}_{ij}=\text{rot}_{X}(\theta(t))\m{b}_{ij},
\end{equation*}
where $\text{rot}_{X}(\cdot)$ is a rotation of a sinusoid angle $\theta(t)=\theta_0\sin(2\pi ft)$ about an arbitrary direction $X\in \mb{R}^3$ (e.g., the principle axes of the coordinate system). The magnitude of $\theta_0$ is up to $10\deg$ and the frequency $f$ is ranged from $1$ to $25$ Hz. Simulation results are provided in Fig. \ref{fig:sim_noise}. It is observed that the estimates of poses of the agents converge to a neighborhood of the actual poses asymptotically.

\section{Conclusion}\label{sec:Conclusion}
In this paper, a network pose localization scheme was proposed for \textit{twin-leader-follower} networks by using direction measurements in $\mb{R}^3$. 
In particular, an orientation localization law in $SO(3)$ and a position localization protocol were presented. We showed that the actual orientations and positions of all follower agents can be estimated almost globally and asymptotically. An extension of this work to systems with more general graph topologies is left as future work.

\begin{ack}                               
{The work of  Q. V. Tran and H.-S. Ahn is supported by the National Research Foundation (NRF) of Korea under the grant NRF2017R1A2B3007034, and the work of B. D. O. Anderson is supported by the Australian Research Council (ARC) under grants DP-160104500 and DP190100887.}
\end{ack}

\bibliographystyle{unsrt}        
\bibliography{quoc2018,quoc2019}           

\begin{thebibliography}{10}

\bibitem{Oh2015survey}
K.-K. Oh, M.-C. Park, and H.-S. Ahn.
\newblock A survey of multi-agent formation control.
\newblock {\em Automatica}, 53:424--440, 2015.

\bibitem{Zhao2018csm}
S.~Zhao and D.~Zelazo.
\newblock Bearing rigidity theory and its applications for control and
  estimation of network systems: Life beyond distance rigidity.
\newblock {\em IEEE Control Syst. Mag.}, 39(2):66--83, 2019.

\bibitem{Aspnes2006}
J.~Aspnes, T.~Ren, D.~K. Goldenberg, A.~S. Morse, W.~Whiteley, Y.~R. Yang,
  B.~D.~O. Anderson, and P.~N. Belhumeur.
\newblock A theory of network localization.
\newblock {\em IEEE Transactions on Mobile Computing}, 40(5):1663--1678, 2006.

\bibitem{Mao2007}
G.~Mao, B.~Fidan, and B.~D.~O. Anderson.
\newblock Wireless sensor network localization techniques.
\newblock {\em Computer Networks}, 51:2529--2553, 2007.

\bibitem{Ma2004}
Y.~Ma, S.~Soatto, J.~Kosecka, and S.~Sastry.
\newblock {\em An Invitation to 3D Vision}.
\newblock Springer, New York, 2004.

\bibitem{Piovan2013}
G.~Piovan, I.~Shames, B.~Fidan, F.~Bullo, and B.~D.O. Anderson.
\newblock On frame and orientation localization for relative sensing networks.
\newblock {\em Automatica}, 49:206--213, 2013.

\bibitem{Tron2014tac}
R.~Tron and R.~Vidal.
\newblock Distributed $3$-{D} localization of camera sensor networks from
  $2$-{D} image measurements.
\newblock {\em IEEE Trans. Autom. Control}, 59(12):3325--3340, 2014.

\bibitem{Montijano2016tro}
E.~Montijano, E.~Cristofalo, D.~Zhou, M.~Schwager, and C.~Sag\"u\'es.
\newblock Vision-based distributed formation control without an external
  positioning system.
\newblock {\em IEEE Trans. Rob.}, 32(2):339--351, 2016.

\bibitem{Zhu2008}
Y.~Zhu, D.~Huang, and A.~Jiang.
\newblock Network localization using angle of arrival.
\newblock In {\em Proc. IEEE Confer. Electro/Information Technology}, pages
  205--210, 2008.

\bibitem{Rong2006}
P.~Rong and M.L. Sitichiu.
\newblock Angle of arrival localization for wireless sensor networks.
\newblock In {\em Proc. 3rd IEEE Confer. Sensor and Ad hoc Communications and
  Networks}, pages 374--382, 2006.

\bibitem{Bishop2011ifac}
A.~N. Bishop and I.~Shames.
\newblock Noisy network localization via optimal measurement refinement part 1:
  {Bearing-only orientation registration and localization}.
\newblock {\em IFAC Proceedings Volumes}, 44(1):8842--8847, 2011.

\bibitem{Thunberg2017cdc}
J.~Thunberg, F.~Bernard, and J.~Gon\c{c}alves.
\newblock Distributed synchronization of {E}uclidean transformations with
  guaranteed convergence.
\newblock In {\em Proc. 56th IEEE Confer. Decision Control (CDC)}, pages
  3757--3762, 2017.

\bibitem{Quoc2018cdc}
Q.~V. Tran, H.-S. Ahn, and B.~D.~O. Anderson.
\newblock Distributed orientation localization of multi-agent systems in
  $3$-dimensional space with direction-only measurements.
\newblock In {\em {P}roc. 57th IEEE Confer. Decision Control (CDC)}, pages
  2883--2889, 2018.

\bibitem{Quoc2018tcns}
Q.~V. Tran, M.~H. Trinh, D.~Zelazo, D.~Mukherjee, and H.-S. Ahn.
\newblock Finite-time bearing-only formation control via distributed global
  orientation estimation.
\newblock {\em IEEE Trans. Control Network Syst.}, 2(6):702--712, 2019.

\bibitem{Zhao2016aut}
S.~Zhao and D.~Zelazo.
\newblock Localizability and distributed protocols for bearing-based network
  localization in arbitrary dimensions.
\newblock {\em Automatica}, 69:334--341, 2016.

\bibitem{Lee2016auto}
B.-H. Lee and H.-S. Ahn.
\newblock Distributed control via global orientation estimation.
\newblock {\em Automatica}, 73:125--129, 2016.

\bibitem{Oh2014tac}
K.-K. Oh and H.-S. Ahn.
\newblock Formation control and network localization via orientation alignment.
\newblock {\em IEEE Trans. Autom. Control}, 59(2):540--545, 2014.

\bibitem{Wang2018tsp}
B.~Wang and Y.-P. Tian.
\newblock Distributed network localization: Accurate estimation with noisy
  measurement and communication information.
\newblock {\em IEEE Transactions on Signal Processing}, 66(22):5927--5940,
  2018.

\bibitem{Bullo2005spr}
F.~Bullo and A.~Lewis.
\newblock {\em Geometric Control of Mechanical Systems}.
\newblock Texts in Applied Mathematics. Springer-Verlag, New York, 2005.

\bibitem{Mahony2008tac}
R.~Mahony, T.~Hamel, and J.-M. Pflimlin.
\newblock Nonlinear complementary filters on the special orthogonal group.
\newblock {\em IEEE Trans. Autom. Control}, 53(5):1203--1218, 2008.

\bibitem{Minh2017tac}
M.~H. Trinh, S.~Zhao, Z.~Sun, D.~Zelazo, B.~D.~O. Anderson, and H.-S. Ahn.
\newblock Bearing-based formation control of a group of agents with
  leader-first follower structure.
\newblock {\em IEEE Trans. Autom. Control}, 64(2):598--613, 2019.

\bibitem{Wahba1965}
G.~Wahba.
\newblock Problem 65–1: A least squares estimate of satellite attitude.
\newblock {\em SIAM Review}, 7(5):409, 1965.

\bibitem{Gower2004}
J.~C. Gower and G.~B. Dijksterhuis.
\newblock {\em Procrustes Problems}.
\newblock Oxford University Press, 2004.

\bibitem{TLee2015tac}
T.~Lee.
\newblock Global exponential attitude tracking controls on {SO(3)}.
\newblock {\em IEEE Trans. Autom. Control}, 60(10):2837--2842, 2015.

\bibitem{Izadi2016auto}
M.~Izadi and A.~K. Sanyal.
\newblock Rigid body pose estimation based on the {L}agrange-d'{A}lembert
  principle.
\newblock {\em Automatica}, 71:78--88, 2016.

\bibitem{Irving2004}
R.~S. Irving.
\newblock {\em Integers, Polynomials, and Rings}.
\newblock Springer-Verlag, New York, 2004.

\bibitem{Caron2005}
R.~Caron and T.~Traynor.
\newblock The zero set of a polynomial.
\newblock {\em University of Windsor}, 2005.
\newblock Technical Report.

\bibitem{Angeli2011tac}
D.~Angeli and L.~Praly.
\newblock Stability robustness in the presence of exponentially unstable
  isolated equilibria.
\newblock {\em IEEE Trans. Autom. Control}, 56(7):1582--1592, 2011.

\bibitem{Angeli2004tac}
D.~Angeli.
\newblock An almost global notion of input to state stability.
\newblock {\em IEEE Trans. Autom. Control}, 6(49):866--874, 2004.

\bibitem{Khalil2002}
H.~K. Khalil.
\newblock {\em Nonlinear Systems}.
\newblock Prentice Hall, third edition, 2002.

\end{thebibliography}



\appendix
\section{Proof of Lemma \ref{lm:left_der_error_function}}\label{app:left_der_error_proof}
Considering each term under the sum in \eqref{eq:orient_cost_function}, we have
\begin{align*}
&\m{D}_{\hat{\m{R}}_i}\Phi_{ij}(\hat{\m{R}}_i,\m{R}_i)\cdot\hat{\m{R}}_i\eta_i^\wedge=-\frac{d}{d\varepsilon}\Big\rvert_{\varepsilon=0}\Phi_{ij}(\hat{\m{R}}_iexp(\varepsilon\eta_i^\wedge),\m{R}_i)\\
&\quad=-tr(\hat{\m{R}}_i\eta^\wedge\m{R}_i^\top\m{b}_{ij}\m{b}_{ij}^\top)=-tr(\eta^\wedge\m{R}_i^\top\m{b}_{ij}\m{b}_{ij}^\top\hat{\m{R}}_i)\\
&\quad =\eta_i^\top(\m{b}_{ij}^i\m{b}_{ij}^\top\hat{\m{R}}_i-\hat{\m{R}}_i^\top\m{b}_{ij}(\m{b}_{ij}^i)^\top)^\vee\\
&\quad = \eta_i^\top(\hat{\m{R}}_i^\top\m{b}_{ij}\times \m{b}_{ij}^i),
\end{align*}
where the third and fourth equalities are derived using \eqref{eq:trace_property_2} and \eqref{eq:cross_prod_3}, respectively. This shows $\m{e}_{ij}$ is a (left-trivialized) derivative of $\Phi_{ij}(\hat{\m{R}}_i)$ with respect to $\hat{\m{R}}_i$. Similarly, we can show that 
\begin{equation*}
\m{D}_{{\m{R}}_i}\Phi_{ij}\cdot{\m{R}}_i\zeta_i^\wedge=-\zeta_i^\top(\hat{\m{R}}_i^\top\m{b}_{ij}\times \m{b}_{ij}^i),
\end{equation*}
and the proof is completed. 
\eprf
\section{Proof of Lemma \ref{lm:time_der_error_function}}\label{app:time_der_error_func_proof}
First, using \eqref{eq:cross_prod_3} we rewrite $\m{e}_i$ as
\begin{align*}
\m{e}_i^\wedge &=\sum_{j\in \mc{N}_i}k_{ij}(\hat{\m{R}}_i^\top\m{b}_{ij}\times \m{b}_{ij}^i)^\wedge=\sum_{j\in \mc{N}_i}k_{ij}(\hat{\m{R}}_i^\top\m{b}_{ij}\times \m{R}_i^\top\m{b}_{ij})^\wedge\\
&=\sum_{j\in \mc{N}_i}k_{ij}(\m{R}_i^\top\m{b}_{ij}\m{b}_{ij}^\top\hat{\m{R}}_i-\hat{\m{R}}_i^\top\m{b}_{ij}\m{b}_{ij}^\top\m{R}_i). \numberthis \label{eq:e_i_hat}
\end{align*}
One has
\begin{align*}
\dot{\m{e}}_i&=\sum_{j\in \mc{N}_i}k_{ij}(-\Omega_i^\wedge\hat{\m{R}}_i^\top\m{b}_{ij}\times \m{R}_i^\top\m{b}_{ij}-\hat{\m{R}}_i^\top\m{b}_{ij}\times (\omega_i^i)^\wedge\m{R}_i^\top\m{b}_{ij})
\\
&=\sum_{j\in \mc{N}_i}k_{ij}\big((-\omega_i^i+\tilde{\Omega}_i)^\wedge\hat{\m{R}}_i^\top\m{b}_{ij}\times \m{R}_i^\top\m{b}_{ij}\\
&\qquad+(\omega_i^i)^\wedge\m{R}_i^\top\m{b}_{ij}\times\hat{\m{R}}_i^\top\m{b}_{ij}\big)
\\
&=\sum_{j\in \mc{N}_i}k_{ij}\big(\tilde{\Omega}_i^\wedge\hat{\m{R}}_i^\top\m{b}_{ij}\times \m{R}_i^\top\m{b}_{ij}-(\omega_i^i)^\wedge\hat{\m{R}}_i^\top\m{b}_{ij}\times \m{R}_i^\top\m{b}_{ij}\\
&\qquad+(\omega_i^i)^\wedge\m{R}_i^\top\m{b}_{ij}\times\hat{\m{R}}_i^\top\m{b}_{ij}\big)\\
&=\sum_{j\in \mc{N}_i}k_{ij}\big(\tilde{\Omega}_i^\wedge\hat{\m{R}}_i^\top\m{b}_{ij}\times \m{R}_i^\top\m{b}_{ij}+\omega_i^i\times (\hat{\m{R}}_i^\top\m{b}_{ij}\times \m{R}_i^\top\m{b}_{ij})\big)\\
&=\sum_{j\in \mc{N}_i}k_{ij}(\tilde{\Omega}_i^\wedge\hat{\m{R}}_i^\top\m{b}_{ij}\times \m{R}_i^\top\m{b}_{ij})-\omega_i^i\times\m{e}_i,
\end{align*} 
where the second and the forth equalities follow from $\Omega_i=\omega_i^i-\tilde{\Omega}_i$ and the Jacobi identity, i.e., $\m{a}\times(\m{b}\times \m{c})+\m{b}\times(\m{c}\times \m{a})+\m{c}\times(\m{a}\times \m{b})=0$, respectively. 
Thus, $||\dot{\m{e}}_{i}||\leq \sum_{j\in \mc{N}_i}k_{ij}||\tilde{\Omega}_i||+\bar{\omega}_i||\m{e}_i||$ which shows (i).

(ii) is followed directly from Lemma \ref{lm:left_der_error_function}.

We show (iii) as follows. It is shown in Lemma \ref{lm:critical_points} that the error function has the form
\begin{equation*}\Phi_i=\text{tr}(\m{G}(\m{I}_3-\m{P})),
\end{equation*}
where $\m{P}\triangleq\m{U}^\top\tilde{\m{Q}}_i\m{U}$ and $\m{G}=\text{diag}\{\lambda(\m{K}_i)\}$. Moreover, from \eqref{eq:e_i_hat} we have
\begin{align*}
\m{e}_i^\wedge &=
\sum_{j\in \mc{N}_i}k_{ij}(\m{R}_i^\top\m{b}_{ij}\m{b}_{ij}^\top\hat{\m{R}}_i-\hat{\m{R}}_i^\top\m{b}_{ij}\m{b}_{ij}^\top\m{R}_i)\\
&=\sum_{j\in \mc{N}_i}(k_{ij}\m{R}_i^\top\m{b}_{ij}\m{b}_{ij}^\top\hat{\m{R}}_i-\hat{\m{R}}_i^\top k_{ij}\m{b}_{ij}\m{b}_{ij}^\top\m{R}_i)\\
&=\m{R}_i^\top\sum_{j\in \mc{N}_i}(k_{ij}\m{b}_{ij}\m{b}_{ij}^\top\hat{\m{R}}_i\m{R}_i^\top-\m{R}_i\hat{\m{R}}_i^\top k_{ij}\m{b}_{ij}\m{b}_{ij}^\top)\m{R}_i\\
&=\m{R}_i^\top( \m{K}_i\tilde{\m{Q}}_i-\tilde{\m{Q}}_i^\top\m{K}_i)\m{R}_i\\
&=\m{R}_i^\top(\m{UG}\m{U}^\top\tilde{\m{Q}}_i-\tilde{\m{Q}}_i^\top\m{UG}\m{U}^\top)\m{R}_i\\
&=\m{R}_i^\top\m{U}(\m{G}\m{U}^\top\tilde{\m{Q}}_i\m{U}-\m{U}^\top\tilde{\m{Q}}_i^\top\m{UG})\m{U}^\top\m{R}_i\\
&=\m{R}_i^\top\m{U}(\m{G}\m{P}-\m{P}^\top\m{G})\m{U}^\top\m{R}_i
\end{align*}
Then $\Phi_i$ is bounded by the square norm of, $||\m{e}_i||=||(\m{G}\m{P}-\m{P}^\top\m{G})^\vee||$, i.e.,
\begin{equation}
\sigma_i||\m{e}_i||^2\leq\Phi_i(\hat{\m{R}}_i, \m{R}_i)\le \gamma_i||\m{e}_i||^2
\end{equation}
where $\sigma_i(\lambda(\m{K}_i)),\gamma_i(\lambda(\m{K}_i))>0$ and the upper bound holds when 
${\Phi}_i<2\min \{\lambda_1+\lambda_2,\lambda_1+\lambda_3,\lambda_2+\lambda_3\},~\lambda_k=\lambda(\m{K}_i),k=1,2,3$ \cite[Prop. 1]{TLee2015tac}.
\eprf
\section{Proof of Theorem \ref{thm:orient_agent_3}}\label{app:orient_agent_3_proof}
We first show that the desired equilibrium point of \eqref{eq:orient_est_law_3} is asymptotically stable. In the second step, we show that the undesired equilibria of \eqref{eq:orient_est_law_3} are unstable.
\subsection*{Step 1: Asymptotic Stability of the Desired Equilibrium}
Consider the Lyapunov function 
\begin{equation}\label{eq:cost_funct_3}
V_3=\frac{1}{2}\tilde{\Omega}_3^2+\Phi_3(\hat{\m{R}}_3, \m{R}_3)-k_V\tilde{\Omega}_3\cdot\m{e}_3,
\end{equation} for a positive constant $k_V$,
which is continuously differentiable and radially unbounded. Following Lemma \ref{lm:time_der_error_function}(iii), we can show that
\begin{equation*}
V_3\geq\frac{1}{2}\m{z}_3^\top\left[\begin{matrix}
1 &-k_V\\
-k_V &2\sigma_3
\end{matrix}\right]
\m{z}_3,
\end{equation*}
where $\m{z}_3^\top=\big[||\tilde{\Omega}_3||,||\m{e}_3||\big]$.
It follows that $V_3\geq0$ if and only if $2\sigma_3-k_V^2>0\leftrightarrow k_V<\sqrt{2\sigma_3}$. 

The time derivative of $V_3$ along the trajectory of \eqref{eq:orient_est_law_3} is given as
\begin{align*}
\dot{V}_3&=-k_\omega\tilde{\Omega}_3^\top\tilde{\Omega}_3-k_V\dot{\tilde{\Omega}}_3\cdot\m{e}_3-k_V\tilde{\Omega}_3\cdot\dot{\m{e}}_3\\
&=-k_\omega\tilde{\Omega}_3^\top\tilde{\Omega}_3-k_V(-k_\omega\tilde{\Omega}_3+\m{e}_3)\cdot\m{e}_3-k_V\tilde{\Omega}_3\cdot\dot{\m{e}}_3\\
&\leq -k_\omega\tilde{\Omega}_3^\top\tilde{\Omega}_3+k_Vk_\omega\tilde{\Omega}_3\cdot \m{e}_3-k_V\m{e}_3^2\\
&\qquad+k_V\tilde{\Omega}_3^\top\Big[\sum_{j\in \mc{N}_3}k_{3j}||\tilde{\Omega}_3||+\bar{\omega}_3||\m{e}_3||\Big]
\\
&\leq -\big[k_\omega-k_V(\sum_{j\in \mc{N}_3}k_{3j})\big]\tilde{\Omega}_3^2+k_V(k_\omega+\bar{\omega}_3)\tilde{\Omega}_3\cdot \m{e}_3-k_V\m{e}_3^2
\\
&\leq -\frac{1}{2}\left[\begin{matrix}
||\m{e}_3|| &||\tilde{\Omega}_3||
\end{matrix}\right]\m{M}_3
\left[\begin{matrix}
||\m{e}_3||\\
||\tilde{\Omega}_3||
\end{matrix}\right],
\end{align*}
where, 
\begin{equation*}
\m{M}_3=\left[\begin{matrix}
2k_V &-k_V(k_\omega+\bar{\omega}_3)\\
-k_V(k_\omega+\bar{\omega}_3) &2(k_\omega-k_V(\sum_{j\in \mc{N}_3}k_{3j}))
\end{matrix}\right].
\end{equation*}
It follows that $\dot{V}_3$ is negative definite if and only if $4k_V(k_\omega-k_V(\sum_{j\in \mc{N}_i}k_{ij}))-k_V^2(k_\omega+\bar{\omega}_3)^2>0\leftrightarrow k_V<\frac{4k_\omega}{(k_\omega+\bar{\omega}_3)^2+4\sum_{j\in \mc{N}_3}k_{3j}}$. Thus, if we choose $k_V$ such that 
\begin{equation}
k_V< \min\left\{\sqrt{2\sigma_3},\frac{4k_\omega}{(k_\omega+\bar{\omega}_3)^2+4\sum_{j\in \mc{N}_3}k_{3j}} \right\}
\end{equation}
then $V_3$ is positive definite and $\dot{V}_3(t)$ is negative definite.
This bounds $V_3(t)\leq V_3(0)$ and consequently $\tilde{\Omega}_3$ is bounded. A direct calculation of $\ddot{V}_3$ shows that $\ddot{V}_3$ is bounded due to boundednesses of $\dot{\tilde{\Omega}}_3$ and $\dot{\m{e}}_3$ (Lemma \ref{lm:time_der_error_function}(i)). As a result, $\tilde{\Omega}_3(t)\rightarrow \m{0},\m{e}_3(t)\rightarrow\m{0}$ as $t\rightarrow\infty$ according to Barbalat's lemma. Consequently, the equilibrium point of \eqref{eq:orient_est_law_3} satisfies $\tilde{\Omega}_3=\m{0}$ and $\tilde{\m{Q}}_3$ is a critical point of $\Phi_3(\tilde{\m{Q}}_3)$ (given in Lemma \ref{lm:critical_points}) where its derivative vanishes (since $\dot{\Phi}_3(\tilde{\m{Q}}_3)=-\tilde{\Omega}_3\cdot \m{e}_3\rightarrow 0$ by Lemma \ref{lm:time_der_error_function}); hence (i) is proved.

It can be shown that the Hessian of $\Phi_3(\tilde{\m{Q}}_3)=\text{tr}(\m{G}(\m{I}_3-\m{U}^\top\tilde{\m{Q}}_3\m{U}))$ at the desired equilibrium $\tilde{\m{Q}}_3=\m{I}_3$ is positive definite \cite[Prop. 11.31]{Bullo2005spr}. Thus, $\hat{\m{R}}_3=\m{R}_3$ is the global minimum of $\Phi_3(\hat{\m{R}}_3)$ (as it will be shown below the other points are either maximum or saddle points). Consequently, $(\hat{\m{R}}_3=\m{R}_3,\tilde{\Omega}_3=\m{0})$ is (locally) asymptotically stable.

\subsection*{Step 2: Instability of Three Undesired Equilibria}
We now show that three undesired equilibria are unstable in what follows. 
\begin{Lemma}\label{lm:hyper_points}
Consider three undesired equilibrium points of \eqref{eq:orient_est_law_3}, i.e., $\{(\tilde{\m{Q}}_3,\tilde{\Omega}_3)~| ~\tilde{\m{Q}}_3\in \{\m{U}\m{D}_1\m{U}^\top,\m{U}\m{D}_2\m{U}^\top,$ $\m{U}\m{D}_3\m{U}^\top\},\tilde{\Omega}_3=\m{0}\}$, where $\m{D}_i,i=1,2,3$ and $\m{U}$ are defined in Lemma \ref{lm:critical_points}. Then, we have $\tilde{\m{Q}}_3$ is either a global maximum or a saddle point of $\Phi_3(\tilde{\m{Q}}_3)=\text{tr}(\m{G}(\m{I}_3-\m{U}^\top\tilde{\m{Q}}_3\m{U}))$.    
\end{Lemma}
\begin{pf}
Consider the variation of $\Phi_3(\m{P}_3)=\text{tr}(\m{G}(\m{I}_3-\m{P}_3))$, where $\m{P}_3:=\m{U}^\top\tilde{\m{Q}}_3\m{U}$, with respect to $\delta\m{P}_3=\m{P}_3\eta^\wedge,\eta^\wedge\in \mathfrak{so}(3)$ as follows
\begin{align*}
\partial_{\m{P}_3}\Phi_3&=\text{tr}(\m{G}(-\m{P}_3\eta^\wedge))=-\text{tr}(\eta^\wedge\m{G}\m{P}_3)\\
&=\eta^\top(\m{G}\m{P}_3-\m{P}_3^\top\m{G})^\vee.
\end{align*}
The second variation of $\Phi_3(\m{P}_3)$ with respect to $\delta\m{P}_3=\m{P}_3\eta^\wedge$ is given as
\begin{align*}
\partial^2_{\m{P}_3}\Phi_3&=-\text{tr}(\delta\eta^\wedge\m{G}\m{P}_3)-\text{tr}(\eta^\wedge\m{G}\m{P}_3\eta^\wedge)\\
&=-\text{tr}(\m{G}\m{P}_3(\eta^\wedge)^2),
\end{align*}
where the second equality follows from the facts that $\delta\eta^\wedge$ is skew-symmetric and $\m{G}\m{P}_3$ is symmetric at equilibrium of points of $\Phi_3$ (due to the condition $\partial_{\m{P}_3}\Phi_3=0$). Direct calculations of $\text{tr}(\m{G}\m{P}_3(\eta^\wedge)^2)$ show that
\begin{align*}
-\text{tr}(\m{G}\m{P}_3(\eta^\wedge)^2)&=\eta^\top\big(\text{tr}(\m{G}\m{P}_3)\m{I}_3-\m{G}\m{P}_3\big)\eta.
\end{align*}
which is the Hessian of $\Phi_3$ evaluated at the critical points of $\Phi_3$.
\\
Consider $\m{P}_3^\prime=\m{U}^\top\m{U}\m{D}_1\m{U}^\top\m{U}=\text{diag}\{1,-1,-1\}$ (the following arguments apply for the other points similarly). Then, 
\begin{align*}
\text{tr}(\m{G}\m{P}_3^\prime)\m{I}_3-\m{G}\m{P}_3\prime=\text{diag}\{-\lambda_2-\lambda_3,\lambda_1-\lambda_3,\lambda_1-\lambda_2\},
\end{align*}
which shows that $\m{P}_3^\prime$ is either a global maximum or a saddle point of $\Phi_3$ depending on the distinct eigenvalues of $\m{K}_3$ defined below Eq. \eqref{eq:error_function_rewrite}.
\eprf
\end{pf}
Consider the error function evaluated at the first undesired equilibrium 
\begin{align*}
\Phi_3(\m{U}\m{D}_1\m{U}^\top)&=\text{tr}(\m{G}(\m{I}_3-\m{D}_1))\\
&=\text{tr(diag}\{\lambda_1,\lambda_2,\lambda_3\}\text{diag}\{0,2,2\})\\
&=2(\lambda_2+\lambda_3).
\end{align*}
Consider the Lyapunov function $U_3=2(\lambda_2+\lambda_3)-V_3$. It holds that  $U_3(\tilde{\m{Q}}_3=\m{U}\m{D}_1\m{U}^\top, \tilde{\Omega}_3=\m{0})=0$. If $\Omega_3$ is sufficiently small we can choose $\hat{\m{R}}_3$ arbitrary close to $\m{R}_i\tilde{\m{Q}}_3$ such that $U_3>0$ (due to Lemma \ref{lm:hyper_points}). Moreover, $\dot{U}_3=-\dot{V}_3>0$. It follows from Chetaev's theorem \cite[Thm. 4.3]{Khalil2002} that $(\tilde{\m{Q}}_3=\m{U}\m{D}_1\m{U}^\top, \tilde{\Omega}_3=\m{0})$ is unstable. Using a similar argument we also conclude that the other undesired equilibrium points are unstable.

Consequently, the desired equilibrium is globally asymptotically stable except a set of measure zero in $SO(3)$ which contains the stable manifolds of the undesired equilibrium points. 
This shows (ii).\eprf

\section{Proof of Lemma \ref{lm:ultimate_boundedness}}\label{app:ultimate_boundedness}
\subsection{Input-To-State Stability}
Consider the Lyapunov function which is similar to \eqref{eq:cost_funct_3} as 
\begin{equation}\label{eq:cost_function_k}
V_k=\frac{1}{2}\tilde{\Omega}_k^2+\Phi_k-k_V\tilde{\Omega}_k\cdot\m{e}_k.
\end{equation}
Then, from Lemma \ref{lm:time_der_error_function}(iii), in $\mc{S}_k$ one has
\begin{equation}\label{eq:bounded_Vk}
1/2\m{z}_k^\top\m{A}_k\m{z}_k\leq V_k \leq 1/2\m{z}_k^\top\m{B}_k\m{z}_k,
\end{equation}
where $\m{z}_k^\top=\big[||\tilde{\Omega}_k||,||\m{e}_k||\big]$ and 
\begin{equation*}
\m{A}_k=\left[\begin{matrix}
1 &-k_V\\
-k_V &2\sigma_k
\end{matrix}\right],
\m{B}_k=\left[\begin{matrix}
1 &-k_V\\
-k_V &2\gamma_k
\end{matrix}\right].
\end{equation*}
Similar to the proof of Theorem \ref{thm:orient_agent_3}, the time derivative of $V_k$ along the trajectory of \eqref{eq:orient_est_law_k} satisfies
\begin{equation}\label{eq:dot_Vk}
\dot{V}_k\leq -\frac{1}{2}\m{z}_k^\top\m{C}_k\m{z}_k+d||\m{h}_k||,
\end{equation}
where $d=\sup\nolimits_t(||\tilde{\Omega}_k-k_V\m{e}_k||)$ and
\begin{equation*}
\m{C}_k=\left[\begin{matrix}
2(k_\omega-k_V(\sum_{j\in \mc{N}_k}k_{kj})) &-k_V(k_\omega+\bar{\omega}_k)\\
-k_V(k_\omega+\bar{\omega}_k) &2k_V
\end{matrix}\right]\footnotemark{}.
\end{equation*}\footnotetext{The scalar $k_{kj}>0,j\in \mc{N}_k$ (defined in \eqref{eq:orient_cost_function}) associates with the edge ($k,j$) where the subscript $\cdot_{kj}$ denotes the agent $k$ and its neighbor $j$.}
Thus, if we choose $k_V$ such that 
\begin{equation}
k_V< \min\left\{\sqrt{2\sigma_k},\frac{4k_\omega}{(k_\omega+\bar{\omega}_k)^2+4\sum_{j\in \mc{N}_k}k_{kj}} \right\}
\end{equation}
then all $\m{A}_k,\m{B}_k$ and $\m{C}_k$ are positive definite. Therefore, it follows from \eqref{eq:bounded_Vk} we have that
\begin{equation}\label{eq:ultimate_bounded_Vk}
\dot{V}_k\leq -\frac{\lambda_{\text{min}}(\m{C}_k)}{\lambda_{\text{max}}(\m{B}_k)}V_k+d||\m{h}_k||.
\end{equation}
The bound of $||\tilde{\Omega}_k||$ can be obtained by considering the derivative of ${U}_k=1/2\tilde{\Omega}_k^2+\Phi_k\geq0$ along the trajectory of \eqref{eq:orient_est_law_k} as follows
\begin{equation*}
\dot{U}_k=-k_\omega\tilde{\Omega}_k^2+\tilde{\Omega}_k^\top\m{h}_k(t)\leq-||\tilde{\Omega}_k||(k_\omega||\tilde{\Omega}_k||-||\m{h}_k||).
\end{equation*}
It follows that $\dot{U}_k\leq 0$ when $||\tilde{\Omega}_k||\geq ||\m{h}_k||/k_\omega$, which shows that $||\tilde{\Omega}_k||$ is ultimately bounded.
\\
Thus, it follows from \eqref{eq:ultimate_bounded_Vk} that the system \eqref{eq:orient_est_law_k} fulfils ultimate boundedness according to \cite[Prop. 3]{Angeli2011tac}.  
This shows input-to-state stability of the unforced system \eqref{eq:orient_est_law_k} with respect to $\m{h}_k(t)$ \cite{Angeli2011tac}.
\subsection{Almost Global Convergence}
It follows from \eqref{eq:ultimate_bounded_Vk} that $\dot{V}_k<0$ if 
\begin{equation}\label{eq:pos_invar_set}
V_k>\frac{\lambda_{\text{max}}(\m{B}_k)}{\lambda_{\text{min}}(\m{C}_k)}d||\m{h}_k||=:\epsilon_1.
\end{equation}
Define the sublevel set $\mc{L}_{\epsilon}:=\{(\tilde{\m{Q}}_k,\tilde{\Omega}_k)\in SO(3)\times \mb{R}^3~|~V_k\leq \epsilon\}$. Then, $\mc{L}_{\epsilon_1}$ is a positive invariant set. Since $\m{h}_k(t)$ tends to zero as $t\rightarrow\infty$, the same is true for $V_k$. To guarantee that $\tilde{\m{Q}}_k\in\mc{S}_k=\{\tilde{\m{Q}}_k~|~\Phi_k(\tilde{\m{Q}}_k)<\phi_k\}$ we consider $V_k< \frac{\phi_k\lambda_{\text{min}}(\m{A}_k)}{2\gamma_k}=:\epsilon_2$. Then, following Lemma \ref{lm:time_der_error_function}(iii) and \eqref{eq:bounded_Vk}, one has
\begin{equation*}
\Phi_k\le \gamma_k||\m{e}_k||^2\leq \gamma_k||\m{z}_k||^2\leq \frac{2\gamma_kV_k}{\lambda_{\text{min}}(\m{A}_k)}<\phi_k.
\end{equation*}
Consequently, any trajectory initializes in or enters $\mc{L}_{\epsilon_2}$ will converge to $\mc{L}_{\epsilon_1}$ and eventually reach $(\tilde{\m{Q}}_k\equiv\m{I}_3,\tilde{\Omega}_k\equiv\m{0})$ as $t\rightarrow\infty$. 
We complete the proof by noting that the sublevel set $\mc{L}_{\epsilon_2}$ covers almost all $SO(3)\times \mb{R}^3$ when the positive scalars $k_{kj}, j\in \mc{N}_k$ are sufficiently large. In particular, as $k_{kj}\rightarrow\infty, \forall j\in \mc{N}_k$, we have that $\phi_k\rightarrow \infty$ and hence $\epsilon_2\rightarrow\infty$. When $\phi_k\rightarrow \infty$ or if $k_{kj}$ are selected such that $\lambda_j(\m{K}_k)\rightarrow \lambda_i(\m{K}_k)$, $i,j\in\{1,2,3\},i\neq j$, we have $\Phi_k(\tilde{\m{Q}}_k(0))<\phi_k$ covers almost all $SO(3)$ excluding the undesired critical points, and the initial $\tilde{\Omega}_k(0)$ satisfies 
\begin{align*}
V_k(0) \leq \lambda_{\text{max}}(\m{B}_k)\m{z}_k^2/2&<\epsilon_2\\
\Rightarrow\tilde{\Omega}_k^2(0)&<2\epsilon_2/\lambda_{\text{max}}(\m{B}_k)-\m{e}_k^2(0)\\
&\leq 2\epsilon_2/\lambda_{\text{max}}(\m{B}_k)-\Phi_k(0)/\gamma_k,
\end{align*} 
covers $\mb{R}^3$ when $\epsilon_2\rightarrow\infty$. 
Furthermore, since the other equilibria are unstable (i.e., arbitrary trajectories close to them will diverge), for almost all initializations of $(\tilde{\m{Q}}_k(0),\tilde{\Omega}_k(0))$ the trajectory of the system \eqref{eq:orient_est_law_k} will converge to the desired equilibrium.
\eprf
\section{Proof of Corollary \ref{coroll:robust_to_noise}}\label{app:robust_to_noise}
The proof follows from the input-to-state stability of the system \eqref{eq:orient_est_law_k} w.r.t. input (Lemma \ref{lm:ultimate_boundedness}). It is also noticed from the proof of Lemma \ref{lm:ultimate_boundedness} that the desired equilibrium of the unforced system of \eqref{eq:orient_est_law_k} is locally exponentially stable. Let $\boldsymbol\delta\in \mb{R}^3$ be the augmented error vector introduced by the direction measurement errors in \eqref{eq:orient_est_law_k} i.e.,
\begin{equation}\label{eq:orient_agent_k_with_noise}
\dot{\hat{\m{R}}}_{k}=\hat{\m{R}}_{k}(\omega_k^k-\tilde{\Omega}_{k})^\wedge, ~\dot{\tilde{\Omega}}_{k}=-k_\omega\tilde{\Omega}_{k}+\m{e}_{k}+\m{h}_{k}(t)+\boldsymbol\delta,
\end{equation} 
where $\m{e}_{k}$ and $\m{h}_{k}$ are defined as in \eqref{eq:orient_est_law_k}, and 
\begin{equation*}
\boldsymbol\delta=\sum_{j\in \mc{N}_{k}}k_{kj}(\hat{\m{R}}_k^\top\hat{\m{R}}_j\tilde{\m{b}}_{kj}^j\times \tilde{\m{b}}_{kj}^k-\hat{\m{R}}_k^\top\hat{\m{R}}_j\m{b}_{kj}^j\times \m{b}_{kj}^k),
\end{equation*}
where $\tilde{\m{b}}$ denotes the noisy estimate of the directional vector $\m{b}$.
\\
Following similar arguments in the proof of Lemma \ref{lm:ultimate_boundedness} we can show that if $||\boldsymbol\delta||$ is sufficiently small the trajectory of \eqref{eq:orient_agent_k_with_noise} converges to a neighborhood of the desired equilibrium point as $t\rightarrow\infty$, i.e.,
\begin{equation}\label{eq:robust_invariant_set}
\left\{ ||\m{z}_k||^2\leq\frac{2\lambda_{\text{max}}(\m{B}_k)}{\lambda_{\text{min}}(\m{C}_k)\lambda_{\text{min}}(\m{A}_k)}d||\boldsymbol\delta||\right\},
\end{equation}
which completes the proof.
\eprf
\end{document}